\documentstyle{article}
\input xypic
\input amssym.def
\input amssym
\author{A.A.~Davydov\thanks{Moscow State University, Moscow, Russia}}
\title{Twisting of monoidal structures}

\newtheorem{prop}{Proposition}
\newtheorem{exam}{Example}
\newcommand{\I}{\it id}

\begin{document}
\maketitle
\begin{abstract}
This article is devoted to the investigation of the deformation (twisting) of monoidal structures, such as the associativity constraint of the monoidal category and the monoidal structure of monoidal functor. The sets of  twistings have a (non-abelian) c
ohomological nature. Using this fact the maps from the sets of twistings to some cohomology groups (Hochschild cohomology of K-theory) are constructed. The examples of monoidal categories of bimodules over some algebra, modules and comodules over bialgeb
ra are examined. We specially concentrate on the case of free tensor category.  
\end{abstract}

\tableofcontents

\section{introduction}
According to Tannaka-Krein theory \cite{dm} a group-like object (group, Lie algebra etc.) is in some sence equivalent to its category of representation with natural tensor product operation plus forgetful functor to vector spaces.
Accurately speaking to extract the group (Lie algebra) from its representations category we must regard some extra-structures (monoidal structures) on the category and the functor (section 2).
\newline
Through this point of view the deformation (quantization) of the group (Lie algebra) is the deformation of its representation category and forgetful functor.
In many cases (semisimple groups and Lie algebras) the deformation does not change the category, tensor product on it and the forgetful functor. So the deformation consists in a change of the monoidal structures of the category and the functor.
\newline
In the present paper we shall deal with only this kind of deformations.
\newline
It appears that this deformations are presented by a twistings of monoidal structures by the automorphisms of some monoidal functors (section 3).
\newline
The set of twistings can be identified with the non-abelian cohomology sets of the cosimplicial complex naturally attached to the monoidal functor (section 4).
\newline
These non-abelian cohomologies can be abelianized by means of infinitesimal methods (section 5) or algebraic K-theory (section 6).
\newline
The examples of monoidal categories of bimodules (section 7), modules and comodules over bialgebra (section 8) and free tensor category (sections 9,10) are examined.
\newline
\newline
The author is pleased to thank P.Deligne, A.Levin, A.Reznikov for useful discussions and especially Yu.Bilu for
calling his attention to the paper \cite{sal}.
\newline
The author was partially supported by ISF grant (n M3E300) and grant of RFFR (n 93-01-1542).
The author is very grateful to the Max-Planck-Institut f\"{u}r Mathematik for stimulating atmosphere.

\section{endomorphisms of functors}
Fix a commutative ring $k$.  
All categories and functors will be $k$-linear. It means that all $Hom$-sets are $k$-modules and the functors induces $k$-linear maps between $Hom$-sets.
\newline
For an associative $k$-algebra $R$ we will denote by $R-Mod$ the category of left $R$-modules and by $R-mod$ the subcategory in $R-Mod$ of finitely generated modules. For example if $k$ is a field then $k-Mod$ is the category of all vector spaces over $k
$ and $k-mod$ is its subcategory of finite-dimensional spaces.
\newline
Another definion of $k$-linearity of the category ${\cal A}$ consits in assigning of the action of $k-Mod$ ($k-mod$) on ${\cal A}$, e.g. the categorical pairing
$$k-Mod \times {\cal A}\rightarrow {\cal A},$$
which is associative respectively tensor product in $k-Mod$.
In this case $k$-linearity of the functor $F$ between $k$-linear categories means that $F$ preservs the action of $k-Mod$.
\newline
It is easy to see that the second version of $k$-linearity implies the first and in the case of categories which have infinite (finite) products. The first version implies the existence of the action of the category, which is non-canonically equivalent t
o $k-Mod$ ($k-mod$).
\newline
Here we shall use the second definition of $k$-linearity.
\newline
Suppose that a faithful (injective on the morphisms) \cite{mcl:cwm} $k$-linear functor $F:{\cal A}\rightarrow k-Mod$ have a left (right) $k$-linear adjoint functor $G:k-Mod\rightarrow {\cal A}$. It means that morphisms of functors
$$\alpha :{\I}_{k-Mod}\rightarrow FG,\quad
\beta :GF\rightarrow {\I}_{\cal A},$$
are given and the compositions of endomorphisms
$$\diagram
&F \rto^{F(\alpha )} & FGF \rto^{\beta F} & F
\enddiagram$$
$$\diagram
&G \rto^{\alpha G} & GFG \rto^{G(\beta )} & G
\enddiagram$$
are identical (the definition of the right adjoint can be obtained from this by replacing $F$ by $G$).
\newline
The $k$-linearity of the functor $G$ involves that 
$$G(V) = E\otimes V ,\quad V\in k-Mod$$
for some object $E$ from ${\cal A}$. We will call this object a {\em generator} ({\em cogenerator} in the case when $G$ is a right adjoint to $F$). The (co)generator $E$ is called {\em projective} ({\em injective}) if the morphism $\beta$ ($\alpha$) of t
he corresponding adjunction is surjective (injective).
\begin{prop}\label{ef}
Let $H:{\cal A}\rightarrow {\cal B}$ be a right (left) exact functor between abelian categories and $E$ is a projective generator (injective cogenerator). Then the homomorphism of algebras
$$End(H)\rightarrow C_{End_{\cal B}(H(E))}(H(End_{\cal A}(E))),$$
which sends the endomorphism $\gamma$ to its specialization ${\gamma}_E$, is an isomorphism .
\end{prop}
Proof:
\newline
Suppose that $E$ is a projective generator (the proof in the case of injective cogenerator can be obtained by the inversion of arrows).
\newline
For any object $A$ from ${\cal A}$ consider the commutative diagram
$$\diagram
& H(A) \dto_{{\gamma}_A} & H(E)\otimes F(A) \lto_{H({\beta}_A )} \dto_{{\gamma}_E} \\
& H(A) & H(E)\otimes F(A) \lto_{H({\beta}_A )}  
\enddiagram$$
The surjectivity of $H(\beta )$ implies that the homomorphism $$End(H)\rightarrow End_{\cal B}(E)$$
is injective.
\newline
Conversely, consider the left exact sequence of functors
$$\diagram
& 0 & H \lto & H(E)\otimes F \lto_{H(\beta )} & H(E)\otimes FK \lto
\enddiagram$$
where $K$ is a kernel of $\beta$. The element $r$ of the centralizer $C_{End_{\cal B}(H(E))}(H(End_{\cal A}(E)))$ defines an endomorphisms of two right functors of this sequence and hence an endomorphism of $H$. It is obvious that the specialization of t
his endomorphism coincides with $r$.

\section{monoidal categories}
The ($k$-linear) {\em tensor product} ${\cal G} \boxtimes {\cal
 H}$ of two ($k$-linear) categories ${\cal G}$ and ${\cal H}$ is a category, whose objects are the pairs of objects of ${\cal G}$ and ${\cal H}$
\begin{displaymath}
(X,Y) \qquad X \in {\cal G}, Y \in {\cal H}
\end{displaymath}
and morphisms are tensor product of morphisms from {\cal G} and {\cal H}
\begin{displaymath}
Hom_{{\cal G} \boxtimes {\cal H}}((X_{1},Y_{1}),(X_{2},Y_{2})) = Hom_{{\cal G}}(X_{1},X_{2}) \otimes_{k} 
Hom_{{\cal H}}(Y_{1},Y_{2})
\end{displaymath}
The {\em tensor product} of functors $F:{\cal G}_{1} \rightarrow {\cal H}_{1}$ and $G:{\cal G}_{2} \rightarrow {\cal H}_{2}$ is a functor 
$F \boxtimes G:{\cal G}_{1} \boxtimes {\cal G}_{2} \rightarrow {\cal H}_{1} \boxtimes {\cal H}_{2}$, 
which sends the pair $(X,Y)$ to $(F(X),G(Y))$.     
The {\em monoidal category} is a category {\cal G} \ with a bifunctor 
\begin{displaymath}
\otimes :{\cal G} \boxtimes {\cal G} \longrightarrow {\cal G} \qquad (X,Y) \mapsto X \otimes Y
\end{displaymath} 
which called {\em tensor (or monoidal) product}. This functor must be equiped with a functorial collection of isomorphisms (so-called {\em associativity constraint})
\begin{displaymath}
\varphi_{X,Y,Z} : X \otimes (Y \otimes Z) \rightarrow (X \otimes Y) \otimes Z \qquad \mbox{for any} \quad X,Y,Z \in {\cal G}
\end{displaymath} 
which satisfies to the following {\em pentagon axiom}:
\newline
the diagram
$$\diagram
& X \otimes (Y \otimes (Z \otimes W)) \dto_{I_X \otimes \varphi_{Y,Z,W} } \rto^{ \varphi_{X,Y,Z \otimes W} } 
& (X \otimes Y) \otimes (Z \otimes W) \rto^{ \varphi_{X \otimes Y,Z,W} } & ((X \otimes Y) \otimes Z) \otimes W \\ 
& X \otimes ((Y \otimes Z) \otimes W) \rrto^{ \varphi_{X,Y \otimes Z,W} } &  & (X \otimes (Y \otimes Z)) \otimes W \uto_{ \varphi_{X,Y,Z} \otimes I_W }
\enddiagram$$
is commutative for any objects $X,Y,Z,W \in {\cal G}$.
\newline
Consider two tensor products of objects $X_{1},...,X_{n}$ from ${\cal G}$ with an arbitrary arrangement of the brackets. 
The coherence theorem \cite{mcl:cwm} states that the pentagon axiom implies the existence of a unique isomorphism between them, which is the composition of the associativity constraints.
\newline
An object 1 together with the functorial isomorphisms 
$$\rho_{X}:X \otimes 1 \rightarrow X \quad 
\lambda_{X}:1 \otimes X \rightarrow X$$
in a monoidal category ${\cal G}$ \ is called a {\em unit} if $\lambda_1 = \rho_1$ and the diagrams
$$\diagram
1\otimes (X\otimes Y) \rrto^{\lambda_{X\otimes Y}} \drto_{\varphi_{1,X,Y}} & & X\otimes Y \\
& (1\otimes X)\otimes Y \urto_{\lambda_{X}\otimes I_{Y}} \\
\enddiagram$$
$$\diagram
X\otimes (1\otimes Y) \rrto^{\rho_X\otimes I_{Y}} \drto_{\varphi_{X,1,Y}} & & X\otimes Y \\
& (X\otimes 1)\otimes Y \urto_{I_X\otimes\lambda_Y} \\
\enddiagram$$
$$\diagram
X\otimes (Y\otimes 1) \rrto^{\rho_{X\otimes Y}} \drto_{\varphi_{X,Y,1}} & & X\otimes Y \\
& (X\otimes Y)\otimes 1 \urto_{I_{X}\otimes\rho_{Y}} \\
\enddiagram$$
are commutative for any $X,Y \in {\cal G}$. 
It is easy to see that the unit object is unique up to isomorphism. 
\newline
A {\em monoidal functor} between monoidal categories ${\cal G}$
and ${\cal H}$ is a functor $F : {\cal G} \longrightarrow 
{\cal H}$ , which is equipped with the functorial collection
of isomorphisms (the so-called {\em monoidal structure}) 
\begin{displaymath}
\gamma_{X,Y} : F(X \otimes Y) \rightarrow F(X) \otimes F(Y)
\qquad \mbox{for any} \quad X,Y \in {\cal G}
\end{displaymath}
for which the following diagram is commutative for any objects $X,Y,Z \in {\cal G}$
$$\diagram & F(X\otimes (Y\otimes Z)) \rto^{c_{X,Y\otimes Z}} \dto_{F(\varphi_{X,Y,Z})} & F(X)\otimes F(Y\otimes Z) \rto^{I\otimes c_{Y,Z}} & F(X)\otimes (F(Y)\otimes F(Z)) \dto_{\psi_{F(X),F(Y),F(Z)}} \\
& F((X\otimes Y)\otimes Z) \rto^{c_{X\otimes Y,Z}} &
F(X\otimes Y)\otimes F(Z) \rto^{c_{X,Y}\otimes I} &
(F(X)\otimes F(Y))\otimes F(Z) \\
\enddiagram$$
A morphism $f : F \rightarrow G$ of monoidal functors $F$ and $G$ is called {\em monoidal} if the diagram
$$\diagram & F(X \otimes Y) \rto^{c_{X,Y}} \dto_{f_{X \otimes Y}} & F(X) \otimes F(Y) \dto_{f_{X} \otimes f_{Y}} \\ & G(X \otimes Y) \rto^{d_{X,Y}} & G(X) \otimes G(Y) 
\enddiagram$$
is commutative for any $X,Y \in {\cal G}$. 
\newline
Denote by $Aut^{\otimes}(F)$ the group of monoidal automorphisms of the monoidal functor $F$.
\newline
Monoidal categories ${\cal G}$ and ${\cal H}$ are {\em equivalent} (as monoidal categories) if there are mutually inverse monoidal functors $F:{\cal G}\rightarrow {\cal H}$ and $G:{\cal H}\rightarrow {\cal G}$ with monoidal isomorphisms $F\circ G \simeq 
{\I}_{\cal H}$ and $G\circ F \simeq {\I}_{\cal G}$. It is easy to see that monoidal categories are equivalent if and only if there is a monoidal functor between them which is an equivalence of categories.
\newline
Let us denote by $Aut^{\otimes}({\cal G})$ the group of monoidal autoequivalences of the monoidal category ${\cal G}$ and by $Aut({\cal G},\otimes )$ the group of autoequivalences of the category ${\cal G}$ which preserves the tensor product $\otimes$. 
\newline
The quotient $\psi\circ{\varphi}^{-1}$ of two associativity constraints $\psi$ and $\varphi$ of the functor $\otimes$ is an isomorphism of the functor $\otimes \circ (\otimes \boxtimes {\it id}) = {\it id}^{\otimes 3}$. Thus we have a map 
\begin{displaymath}
Q : MS({\cal G}) \times MS({\cal G}) \longrightarrow Aut({\it id}^{\otimes 3})
\end{displaymath}
from the direct square of the set $MS({\cal G}) = MS({\cal G}, \otimes)$ of all monoidal structures of the category ${\cal G}$ with the tensor product $\otimes$ (that is the set of all associativity constraints of the tensor product $\otimes$) to the gro
up $Aut({\it id}^{\otimes 3})$ of automorphisms of the functor ${\it id}^{\otimes 3}$.
\newline
Analogously, the qoutient $d \circ c^{-1}$ of two monoidal structures $c$ and $d$ of functor $F$ is an isomorphism of the functor $\otimes \circ (F \boxtimes F) = F^{\otimes 2}$. Thus there is a map 
\begin{displaymath}
Q : MS(F) \times MS(F) \longrightarrow Aut(F^{\otimes 2})
\end{displaymath}  
from the direct square of the set $MS(F)$ of all monidal structures of the functor $F$ to the group $Aut(F^{\otimes 2})$.
\newline
It immediately follows from the definitions that these maps satisfy the conditions 
\begin{displaymath}
Q(x,x) = 1 \qquad \mbox{and} \qquad Q(x,y)Q(y,z) = Q(x,z)
\end{displaymath}
The map 
$$Q(\varphi,?) : MS({\cal G}) \longrightarrow Aut({\it id}^{\otimes 3})$$
is injective for any associativity constraint $\varphi$ of the tensor product $\otimes$ of the category ${\cal G}$.
The image of this map consists of those isomorphisms $\alpha \in Aut({\it id}^{\otimes 3})$ for which the diagram 
$$\diagram
& X\otimes (Y\otimes (Z\otimes W)) \dto_{I\otimes \varphi\circ\alpha} \rto^{\varphi\circ\alpha} 
& (X\otimes Y)\otimes (Z\otimes W) \rto^{\varphi\circ\alpha} & ((X\otimes Y)\otimes Z)\otimes W \\ 
& X\otimes ((Y\otimes Z)\otimes W) \rrto^{\varphi\circ\alpha} &  & (X\otimes (Y\otimes Z))\otimes W \uto_{\varphi\circ\alpha\otimes I}
\enddiagram$$
is commutative for any objects $X,Y,Z,W \in {\cal G}$.
\newline
For such $\alpha$ we call the monoidal category ${\cal G}(\alpha ) = ({\cal G},\varphi)(\alpha ) = ({\cal G},\alpha\varphi)$ a {\em twisted form} ({\em twisting} by $\alpha$) of the monoidal category $({\cal G},\varphi)$.
\newline 
Similarly, the map
$$Q(c,?) : MS(F) \longrightarrow Aut(F^{\otimes 2})$$
is injective for any monoidal structure $c$ of the functor $F$ and its image consists of those isomorphisms $\alpha \in Aut(F^{\otimes 2})$ for which the diagram 
$$\diagram 
& F(X\otimes (Y\otimes Z)) \rto^{c\circ\alpha} \dto_{\varphi
} & F(X)\otimes F(Y\otimes Z) \rto^{I\otimes c\circ \alpha} & F(X)\otimes (F(Y)\otimes F(Z)) \dto_{\psi} \\
& F((X\otimes Y)\otimes Z) \rto^{c\circ\alpha} &
F(X\otimes Y)\otimes F(Z) \rto^{c\circ\alpha\otimes I} &
(F(X)\otimes F(Y))\otimes F(Z) \\
\enddiagram$$
is commutative for any objects $X,Y,Z \in {\cal G}$. 
\newline
The pair $(F,\alpha c) = (F,c)(\alpha) = F(\alpha )$ will be called a {\em twisted form} of the monoidal functor $(F,c)$.
\newline 
The notion of the isomorphism of monoidal functors defines an equivalence relation $"\sim "$ on the set $MS(F)$
$$c \sim d \Leftrightarrow (F,c) \simeq (F,d) \qquad c,d \in  MS(F).$$
Denote by $Ms(F)$ the factor of $MS(F)$ by the relation $\sim$.
\newline
The inclusion $Q(c,?)$ induces an equivalence relation on its image (and even on the whole group $Aut(F^{\otimes 2})$): 
\newline 
$\alpha \sim \beta$ iff there is an automorphism $f$ of the functor $F$ such that the diagram
$$\diagram 
& F(X\otimes Y)\rto^{c\circ\alpha} \dto_{f} & F(X) \otimes F(Y) \dto_{f\otimes f} \\ & F(X \otimes Y) \rto^{c\circ\beta} & F(X) \otimes F(Y) 
\enddiagram$$
is commutative for any objects $X,Y \in {\cal G}$.
It is easy to see that this relation on the group $Aut(F^{\otimes 2})$ comes from the action of the group $Aut(F)$.
\newline
In the case of the identity functor the set $Ms({\I}_{\cal G})$ 
has a group structure and coincides with the kernel of the natural map 
$$Aut^{\otimes}({\cal G})\rightarrow Aut({\cal G},\otimes ).$$ 
\newline
In a similar way we can define the factor $Ms({\cal G})$ by the equivalence relation $"\approx "$ on the set $MS({\cal G})$, which corresponds to the equivalence of monoidal categories. Since the definition of the relation $"\approx "$ uses all autoequiv
alences of the category it is often difficult to verify it. It is useful to consider a more weak relation $"\sim "$ which corresponds to the equivalence by means of identity functor \cite{saa:ct}. 
\newline
Namely, $\varphi\sim\psi$ iff the identity functor $\I_{\cal G}$ can be equipped with a structure of monoidal functor between $({\cal G},\otimes ,\varphi )$ and $({\cal G},\otimes ,\psi )$.
\newline
The factor $MS({\cal G})/\sim$ will be denoted by $\overline{Ms}({\cal G})$. 
\newline
The corresponding relation on the image of inclusion $Q(\varphi,?)$ (on the group $Aut({\I}_{\cal G}^{\otimes 3})$) has the form:
\newline
$\alpha\sim\beta$ iff there is an automorphism $c$ of the functor $F^{\otimes 2}$ such that the diagram
$$\diagram 
& F(X\otimes (Y\otimes Z)) \rto^{c} \dto_{\varphi\circ\alpha} & F(X)\otimes F(Y\otimes Z) \rto^{I\otimes c} & F(X)\otimes (F(Y)\otimes F(Z)) \dto_{\varphi\circ\beta} \\
& F((X\otimes Y)\otimes Z) \rto^{c} &
F(X\otimes Y)\otimes F(Z) \rto^{c\otimes I} &
(F(X)\otimes F(Y))\otimes F(Z) \\
\enddiagram$$
is commutative for any objects $X,Y,Z \in {\cal G}$.
\newline
It is evident that this relation on the group $Aut({\I}_{\cal G}^{\otimes 3})$ comes from the action of $Aut({\I}_{\cal G}^{\otimes 2})$ on it.
\newline
The group $Aut({\cal G},\otimes )$ acts on the set $Ms({\cal G})$ in the following way:
$$F(\varphi )_{X,Y,Z} = F(\varphi_{F^{-1}(X),F^{-1}(X),F^{-1}(Z)})$$ 
where $\varphi\in Ms({\cal G})$ and $F\in Aut({\cal G},\otimes )$. The factor of this action coincides with the set $\overline{Ms}({\cal G})$.

\section{endomorphisms of monoidal functor}
The collection of the endomorphisms algebras of tensor powers of a monoidal functor F can be equipped with the structure of cosimplicial complex $E(F)_{*} = End(F^{\otimes *})$.
\newline 
The image of the coface map
$${\partial}_{n+1}^{i} : End(F^{\otimes n}) \rightarrow End(F^{\otimes n+1}) \qquad i = 0,...,n+1$$
of the endomorphism $\alpha \in End(F^{\otimes n})$ has the following specialization on the objects $X_{1},...,X_{n+1}$: 
$${\partial}_{n+1}^{i}(\alpha )_{X_{1},...,X_{n+1}} = \left\{
\begin{array}{l}
{\phi}_{0}(I_{X_{1}} \otimes \alpha_{X_{2},...,X_{n+1}}){\phi}^{-1}_{0}, \quad i = 0 \\
{\phi}_{i}(\alpha_{X_{1},...,X_{i} \otimes X_{i+1},...,X_{n+1}}){\phi}^{-1}_{i}, \quad 1 \leq i \leq n \\
\alpha_{X_{1},...,X_{n}} \otimes I_{X_{n+1}}, \quad i = n+1
\end{array}
\right.,$$
here ${\phi}_{i}$ is the unique isomorphism between $F(...(X_{1} \otimes...) \otimes X_{n+1}) = F^{\otimes n+1}(X_{1},...,X_{n+1})$ and $F(...(X_{1} \otimes...) \otimes (X_{i} \otimes X_{i+1})) \otimes ...) \otimes X_{n+1})$.
The specialization of the image of the coboundary map
$${\sigma}_{n-1}^{i} : End(F^{\otimes n}) \rightarrow End(F^{\otimes n+1}) \qquad i = 0,...,n-1$$ 
is 
$${\sigma}_{n-1}^{i}(\alpha)_{X_{1},...,X_{n-1}} = {\alpha}_{X_{1},...,X_{i},1,X_{i+1},...,X_{n-1}}.$$
We may also define the zero component of this complex as the  endomorphism algebra of the unit object of the category ${\cal H}$ which can be regarded as the endomorphism algebra of the functor 
$$F^{\otimes 0}:k-Mod \longrightarrow {\cal H},\qquad
F^{\otimes 0}(V) = V\otimes 1$$
The coface maps 
$${\partial}_{1}^{i}:End_{\cal H}(1) \rightarrow End(F^{\otimes 1}),\qquad i=0,1$$
has the form
$${\partial}_{1}^{0}(a) = \rho (I_{1}\otimes a){\rho}^{-1}, \qquad {\partial}_{1}^{1}(a) = \lambda (I_{1}\otimes a){\lambda}^{-1};$$
here $\rho$ and $\lambda$ are the structural isomorphisms of the unit object $1$.
\newline  
The invertible elements of this complex of algebras form a cosimplicial complex of the (generally non-commutative) groups $A_{*}(F) = Aut(F^{\otimes *})$.
Nevertheless, for small $n$ the $n$-th cohomology of the complex $A(F)_{*}$ can be defined.
\newline
Let us consider the cosimplicial complex $E_{*}$ of algebras as a functor \cite{gz}
$$E:\delta \longrightarrow k-Alg$$
from the category $\delta$ of finite ordered sets with  nondecreasing maps to the category $k-Alg$ of $k$-algebras. In particular, $E_{n}$  is the value of the functor $E$ on the ordered set of $n+1$ elements $[n]=\{0,1,...n\}$ and $${\partial}_{n}^{i} =
 E({\partial}_{n}^{i}), \qquad
{\sigma}_{n}^{i} = E({\sigma}_{n}^{i}),$$
where
\newline
${\partial}_{n}^{i}$ is the increasing injection which does not take the value $i\in [n]$, 
\newline
${\sigma}_{n}^{i}$ is the nondecreasing surjection which takes twice the value $i\in [n]$.
\newline
We say that the {\em linking coefficient} of two nondecreasing maps $\tau ,\pi :[l]\rightarrow [m]$ is equal to $n$ iff there are decompositions into the nonintersecting unions 
$$Im(\tau) = A_{1}\cup ...\cup A_{s}, \quad A_{1}<...<A_{s}$$
$$Im(\pi) = B_{1}\cup ...\cup B_{t}, \quad  B_{1}<...<B_{t} \qquad s+t=n+1$$
$$\mbox{and}\quad A_{1}\leq B_{1}\leq A_{2}\leq B_{2}...$$
In particular the linking coefficient of ${\partial}_{n}^{i}$ and ${\partial}_{n}^{i+1}$ is equal to
$$1, \ \mbox{if}\ n=1;\qquad
2, \ \mbox{if}\ n=2;\qquad
2n-3, \ \mbox{if}\ n\geq 3$$    
if $i>j-1$ then the linking coefficient of ${\partial}_{n}^{i}$ and ${\partial}_{n}^{j}$ is equal to 
$$1, \ \mbox{if}\ n=2;\qquad 
2, \ \mbox{if}\ n=3;\qquad 
2n-5, \ \mbox{if}\ n\geq 4.$$ 
The functor $E:\delta \longrightarrow k-Alg$ (the cosimplicial complex of algebras $E_{*}$) will be called {\em n-commutative} if 
$$E(\tau)(a)\quad \mbox{commutes with}\quad E(\pi)(b)\qquad \mbox{for every}\quad a,b\in E(l),$$
for any nondecreasing maps $\tau ,\pi :[l]\rightarrow [m]$ whose linking coefficient is less or iqual to $n$, if $n<3$ and $2n-1$, if $n\geq3$.
\newline
For the cosimplicial complex of groups $A_{*}$ let us denote by ${\it z}^{n}(A_{*}) \subset A_{n}$ the set of solutions of the equaition
$$\prod_{i=0}^{\left[\frac{n+1}{2}\right]} {\partial}^{2i}_{n+1}(a) = \prod_{i=\left[\frac{n}{2}\right]}^{0} {\partial}^{2i+1}_{n+1}(a)$$
and by $"\sim "$ the binary relation on the group $A_{n}$:
$$a\sim b \Leftrightarrow \exists c\in A_{n-1}: \qquad a\left(\prod_{i=0}^{\left[\frac{n}{2}\right]} {\partial}^{2i}_{n}(c)\right) = \left(\prod_{i=\left[\frac{n-1}{2}\right]}^{0} {\partial}^{2i+1}_{n}(c)\right)b$$
\begin{prop}
If $A_{*}$ is a n-commutative cosimplicial complex of algebras, then
\newline
1.$A_m$ is a commutative group for $m<n$ 
\newline
2.${\it z}_{m}(A_{\cdot})$ is a subgroup of $A_m$ for $m\leq n$ 
\newline
3.$"\sim "$ is an equivalence relation on the group $A_{m}$ (induced by the action of $A_{m}$ which preservs ${\it z}_{m}(A_{\cdot})$) for $m\leq n+1$
\end{prop} 
Proof:
\newline
1.For $m<n$ and for any $a,b \in A_{m}$ the commutator $[{\partial}_{m+1}^{i}(a),{\partial}_{m+1}^{i+1}(b)]$ is equal to zero. So
$$[a,b] = {\sigma}_{m}^{i} ([{\partial}_{m+1}^{i}(a),{\partial}_{m+1}^{i+1}(b)]) = 0$$
2,3.For $m\leq n$ and for any $a,b \in A_m$
$$
\prod_{i=0}^{\left[\frac{m+1}{2}\right]} {\partial}^{2i}_{m+1}(a) \prod_{i=0}^{\left[\frac{m+1}{2}\right]} {\partial}^{2i}_{m+1}(b) = \prod_{i=0}^{\left[\frac{m+1}{2}\right]} {\partial}^{2i}_{m+1}(ab)
$$
and
$$\prod_{i=\left[\frac{m}{2}\right]}^{0} {\partial}^{2i+1}_{m+1}(a) \prod_{i=\left[\frac{m}{2}\right]}^{0} {\partial}^{2i+1}_{m+1}(a) = \prod_{i=\left[\frac{m}{2}\right]}^{0} {\partial}^{2i+1}_{m+1}(ab)$$
\newline
\newline
It immediately follows from the definition of the coface maps that the complex of automorphisms $A_{*}(F)$ of a monoidal functor $F$ is 1-commutative. 
\newline
The group 
$${\it z}^{1}(F) = {\it z}^{1}(A_{*}(F)) = \{ \alpha\in Aut(F), {\partial}^{1}_{2}(\alpha ) = {\partial}^{0}_{2}(\alpha ){\partial}^{2}_{2}(\alpha )\}$$
is the same as the group of monoidal automorphisms $Aut^{\otimes}(F)$ of the functor $F$.
\newline
The set 
$${\it z}^{2}(F) = \{ \alpha\in Aut(F^{\otimes 2}), {\partial}^{0}_{3}(\alpha ){\partial}^{2}_{3}(\alpha ) =  {\partial}^{3}_{3}(\alpha ){\partial}^{1}_{3}(\alpha )\}$$
coincides with the image of the map $Q(c,?) : MS(F) \longrightarrow Aut(F^{\otimes 2})$ and the equivalence relation $\sim$ on ${\it z}^{2}(F)$ coincides with the relation on the image of $Q(c,?)$ which correspods to the isomorphism of monoidal functors.
 So the second cohomology of monoidal functor is isomorphic to the set of monoidal structures
$${\it h}^{2}(F) = {\it z}^{2}(F)/\sim = Ms(F).$$
It follows from the functoriality of the automorphisms that the complex of automorphisms $A_{*}(\I_{\cal G})$ of the identity functor $\I$ is 2-commutative.
\newline
The group structure on ${\it h}^{2}({\cal G}) = {\it h}^{2}(\I_{\cal G})$ corresponds to the composition of monoidal functors.
\newline
The set 
$${\it z}^{3}({\cal G}) = {\it z}^{3}(\I_{\cal G}) = \{ \alpha\in Aut(\I^{\otimes 3}), {\partial}^{0}_{4}(\alpha ){\partial}^{2}_{4}(\alpha ){\partial}^{4}_{4}(\alpha ) =  {\partial}^{3}_{3}(\alpha ){\partial}^{1}_{2}(\alpha )\}$$
is nothing alse then the image of the map $Q(c,?) : MS({\cal G}) \longrightarrow Aut(\I^{\otimes 3})$.
\newline 
The equivalence relation $\sim$ on ${\it z}^{3}({\cal G})$ coincides with the relation on the image of $Q(c,?)$ which correspods to the equivalence of monoidal categories by means of the identity functor. So the third cohomology of the monoidal category 
(of the identity functor) is isomorphic to the set of monoidal structures
$${\it h}^{3}({\cal G}) = {\it z}^{3}({\cal G})/\sim = Ms({\cal G}).$$

\section{tangent cohomology}
Let $K$ be a $k$-algebra and ${\cal G}$ be a $k$-linear category. 
\newline
The {\em subcategory of $K$-modules} ${\cal G}_K$ in ${\cal G}$ is a category of pairs $(X,\alpha )$ where $X$ is an object of ${\cal G}$ and $\alpha :K\rightarrow End_{\cal G}(X)$ is a homomorphism of $k$-algebras. The morphism from $(X,\alpha )$ to $(Y
,\beta )$ in ${\cal G}_K$ is a morphism $f$ from $X$ to $Y$ in ${\cal G}$ which preserves the $K$-module structure ($f\alpha (c) = \beta (c)f$ for any $c\in K$).
\newline
The forgetful functor 
$${\cal G}_K\rightarrow {\cal G},\qquad (X,\alpha )\mapsto X$$ has the right adjoint
$${\cal G}\rightarrow {\cal G}_K,\qquad X\mapsto (K\otimes_{k} X, i\otimes I_{X})$$
if the category ${\cal G}$ is sufficiently large (e.g. we can multiply any object by the $k$-module $K$). Anyway we have no problem if $K$ is a finitely generated projective $k$-module. 
\newline
A functor $F:{\cal G}\rightarrow {\cal H}$ induces the functor
$$F_K :{\cal G}_{K}\rightarrow {\cal H}_{K},\qquad F(X,\alpha ) = (F(X), F(\alpha ))$$
such that the diagram of functors
$$\diagram
& {\cal G}_K \rto^{F_K} \dto & {\cal H}_K \dto \\
& {\cal G} \rto^F        & {\cal H}
\enddiagram$$
is commutative.
\begin{prop}
Let $F:{\cal G}\rightarrow {\cal H}$ is a right exact monoidal functor between abelian monoidal categories with right biexact tensor products.
\newline
Then for any $k$-algebra $K$ the homomorphism of cosimplicial complexes
$$Z(K)\otimes_k E^{*}(F)\rightarrow E^{*}(F_{K}),\qquad {c\otimes\gamma}_{(X,\alpha)} = \alpha (c)\gamma$$
is an isomorphism (where $Z(K)$ is the center of the algebra $K$).
\end{prop}
Proof:
It is easy to see that the statement of the proposition follows from the following fact:
\newline
for the right exact functor $F:{\cal G}\rightarrow {\cal H}$ the homomorphism of algebras 
$$Z(K)\otimes_k End(F)\rightarrow End(F_{K}),\qquad {c\otimes\gamma}_{(X,\alpha)} = \alpha (c)\gamma$$
is an isomorphism.
\newline
The proof of this fact is analogous to the proof of the proposition (\ref{ef}).
\newline
\newline
Now we apply the functor structure of $h^{*}(F)$ to define the tangent space (tangent cohomology) and the tangent cone.
\newline
The {\em tangent space} to the functor $h:k-Alg\rightarrow Sets$ from the category of $k$-algebras to the category of sets at the point $x\in h(k)$ is a fibre (over this point) $T_{x}(h(k))$ of the map
$$h(A_2)\rightarrow h(k),$$
which is induced by the homomorphism of algebras
$$A_2 = k[\epsilon|\epsilon^2 = 0]\rightarrow k,\qquad \epsilon\mapsto 0.$$
\begin{prop}
The tangent space of the functor $h^{n}(F)$ of cohomologies of monoidal functor $F$ at the unit point (the class of idetity endomorphism) coincides with the cohomology $H^{n}(F)$ of the cochain complex associated with the cosimplicial complex $E^{*}(F)$ 
of abelian groups (the tangent cohomology of $F$).
\end{prop}
Proof\newline
Since $E^{*}(F_{A_2}) = A_2\otimes_{k}E^{*}(F)$, any element of the kernel of the map $A^{n}(F_{A_2})\rightarrow A^{n}(F)$ has the form $1+\epsilon\alpha$ for some $\alpha\in E^{n}(F)$. 
\newline
The direct checking shows that the cocycle condition for $1+\epsilon\alpha$ is the equation $d_{n}(\alpha) = 0$, where $d_{n} = \sum_{i=0}^{n+1}(-1)^{i}{\partial}_{n+1}^{i}$ is a cochain differential. The cocycles $1+\epsilon\alpha$ and $1+\epsilon\beta$
 is equivalent if and only if $\alpha$ and $\beta$ differs by the coboundary.
\newline
\newline
For example, the first cocycle module $Z^{1}(F)$ of the tangent complex coincides with the module of {\em differentiations} of the monoidal functor $F$
$$Diff(F) = \{l\in End(F), l_{X\otimes Y}=l_{X}\otimes I_Y + I_{X}\otimes l_Y \}$$
which is a Lie algebra respectively the commutator in the algebra $End(F)$. The first tangent cohomology $H^{1}(F)$ is a factoralgebra of this Lie algebra.
\newline
The tensor product of endomorphisms of the funtor $F$
$$End(F^{\otimes n})\otimes End(F^{\otimes m})\rightarrow End(F^{\otimes n+m})$$
defines on $E^{*}(F)$ a structure of cosimplicial algebra
$${\partial}^{i}_{n+m+1}(\alpha\otimes\beta) = 
\left\{
\begin{array}{cl}
{\partial}^{i}_{n+1}(\alpha)\otimes\beta, & i\leq n\\
\alpha\otimes{\partial}^{i-n}_{m+1}(\beta), & i>n
\end{array}
\right.$$
Hence the associated cochain complex $(E^{*}(F),d)$ is a differential graded algebra
$$d(\alpha\otimes\beta) = d(\alpha)\otimes\beta + (-1)^{n}\alpha\otimes d(\beta).$$
In particular, the tangent cohomologies $H^{*}(F)$ form a graded algebra. 
\newline
Let us also note that the multiplication is skew-symmetric on the first component $H^{1}(F)$, because
$$\alpha\otimes\beta + \beta\otimes\alpha = d(-\alpha\beta), \qquad\mbox{for} \alpha ,\beta\in Z^{1}(F).$$
The {\em tangent cone} in the tangent space at the point $x\in h(k)$ to the functor $h:k-Alg\rightarrow Sets$ is the image of the map
$$ker(h(A_{\infty})\rightarrow h(k))\rightarrow ker(h(A_2)\rightarrow h(k)),$$
induced by the homomorphism of algebras
$$A_{\infty} = k[[\epsilon ]]\rightarrow  k[\epsilon|\epsilon^2 = 0] = A_2 ,\qquad \epsilon\mapsto\epsilon.$$
In many cases it is sufficient (and more convinient) to consider the {\em first approximation of the tangent cone} which is the image of the map
$$ker(h(A_{3})\rightarrow h(k))\rightarrow ker(h(A_2)\rightarrow h(k)),$$
induced by the homomorphism of algebras
$$A_3 = k[\epsilon|\epsilon^3 = 0]\rightarrow  k[\epsilon|\epsilon^2 = 0] = A_2 ,\qquad \epsilon\mapsto\epsilon.$$

\section{K-theory}
In \cite{qui} Quillen associated to an abelian (exact) category ${\cal A}$ a topological space $BQ{\cal A}$ such that exact functors between categories define the continuous maps between the corresponding spaces and isomorphisms of functors define the ho
motopies between the corresponding maps. In other words, Quillen's space is a 2-functor from the 2-category of abelian (exact) categories (with isomorphisms of functors as 2-morphisms) to the 2-category of topological spaces.
\newline
Waldhausen \cite{wld} proved that the 2-functor $K = \Omega BQ$ is permutable (in some sence) with the product. Namely, he constructed the continuous map $K({\cal A})\wedge K({\cal B})\rightarrow K({\cal C})$ for any biexact functor ${\cal A}\times {\cal
 B}\rightarrow {\cal C}$.
\newline
The homotopy groups $K_{*}({\cal A}) = \pi_{*}({\cal A})$ of the Waldhausen space $K({\cal A})$  are called {\em algebraic K-theory} of the category ${\cal A}$.
\newline
Now we give the definition of Hochschild cohomology of  topological ring spaces.
\newline
The topological space $K$ with a continuous map $\mu:K\wedge K\rightarrow K$ is called a {\em ring space}. A ring space $(K,\mu )$ is {\em associative} if there is a homotopy  between  $\mu (I\wedge\mu)$ and $\mu (\mu\wedge I)$.
\newline
A {\em bimodule} over an associative ring space $(K,\mu )$ is a space $M$ together with the continuous maps 
$$\nu :M\wedge K\rightarrow M, \qquad\upsilon :K\wedge M\rightarrow M$$
and the homotopies
$$\nu (I\wedge\mu )\rightarrow\nu (\nu\wedge I), \quad
\upsilon (I\wedge\nu )\rightarrow\nu (\upsilon\wedge I), \quad \upsilon (I\wedge\upsilon )\rightarrow\upsilon (\mu\wedge I)$$
For example the space of maps to the ring space is a bimodule over this ring space.
\newline
Denote by $[X,Y]$ the set of homotopy classes of continuous maps from $X$ to $Y$.
\newline
The {\em Hochschild complex} of a ring space $K$ with coefficients in a bimodule space $M$ is a semicosimplicial complex of sets
$$C_{*}(K,M), \qquad C_{n}(K,M) = [K^{\wedge n},M]$$
with the coface maps ${\partial}^{i}_{n}:C_{n-1}(K,M)\rightarrow C_{n}(K,M)$ defined as follows
$${\partial}^{i}_{n}(f) = \left\{ \begin{array}{l}
\upsilon (I\wedge f),\qquad i=0\\
f(I\wedge ...\wedge\mu\wedge ...\wedge I),\qquad 1\leq i\leq n\\
\nu(f\wedge\I),\qquad i=n+1
\end{array}\right.$$
If $M$ is a loop space then $C_{\cdot}(K,M)$ is a complex of groups and these groups are abelian if $M$ is a double loop space. In the second case the cohomology of the cochain complex associated with $C_{*}(K,M)$ will be called {\em Hochschild cohomolog
y} ($HH^{*}(K,M)$) of ring space $K$ with coefficient in the bimodule $M$.
\newline
A ring (bimodule) structure on a space induces a graded ring (bimodule) structure on its homotopy groups.
\newline
The natural map $C_{*}(K,M)$ to the Hochschild complex $C_{*}(\pi_{\cdot}(K),\pi_{\cdot}(M))$ of the ring $\pi_{*}(K)$ with coefficients in the bimodule $\pi_{*}(M)$ induces the homomorphism of Hochschild cohomology
$$HH^{*}(K,M) \longrightarrow HH_{*}(\pi_{*}(K),\pi_{*}(M)).$$
\begin{prop}
Let $F:{\cal G}\rightarrow {\cal H}$ be an exact monoidal functor between abelian monoidal categories with biexact tensor products. Then there is a map of semicosimplicial complexes of groups
$$A_{*}(F) = Aut(F^{\otimes *})\longrightarrow C_{*}(K({\cal G}),\Omega K({\cal H}))$$
which defines the maps
$${\it h}^{i}(F)\longrightarrow HH^{i}(K_{*}({\cal G}),K_{*+1}({\cal H}))$$
\end{prop}
Proof:
\newline
The map $Aut(F^{\otimes *})\longrightarrow C_{*}(K({\cal G}),\Omega K({\cal H}))$ sends the automorphism $\alpha \in Aut(F^{\otimes n})$ to the class of the corresponding autohomotopy of the map $K(F^{\otimes n})$. 
\newline
\newline
The zero component of the previous map 
$${\it h}^{i}(F)\longrightarrow HH^{i}(K_{0}({\cal G}),K_{1}({\cal H}))$$
admits a more explicit describrion. 
\newline
For an automorphism $\alpha \in Aut(F^{\otimes n})$ the class $\left[\alpha_{X_{1},...,X_{n}}\right]\in K_{1}({\cal H})$ depends only of the classes $[X_{1}],...,[X_{n}]\in K_{0}({\cal G})$. Indeed, the exact sequence 
$$0\rightarrow Y_{i}\rightarrow X_{i}\rightarrow Z_{i}\rightarrow 0$$
can be extended to the diagram with exact rows
$$\diagram
&0 \rto & X_{1}\otimes...\otimes Y_{i}\otimes...\otimes X_{n} \rto \dto_{\alpha} & X_{1}\otimes...\otimes X_{i}\otimes...\otimes X_{n} \rto \dto_{\alpha} &  X_{1}\otimes...\otimes Z_{i}\otimes...\otimes X_{n} \rto \dto_{\alpha} & 0 \\
&0 \rto & X_{1}\otimes...\otimes Y_{i}\otimes...\otimes X_{n} \rto & X_{1}\otimes...\otimes X_{i}\otimes...\otimes X_{n} \rto &  X_{1}\otimes...\otimes Z_{i}\otimes...\otimes X_{n} \rto & 0 
\enddiagram$$
which is commutative by the functoriality of $\alpha$. Hence
$$\left[\alpha_{X_{1},...,X_{i},...,X_{n}}\right] = \left[\alpha_{X_{1},...,Y_{i},...,X_{n}}\right] + \left[\alpha_{X_{1},...,Z_{i},...,X_{n}}\right]$$
which means that $\left[\alpha_{X_{1},...,X_{n}}\right]$ depends only of the classes $[X_{1}],...,[X_{n}]$ and can be regarded as an element $\overline{\alpha}$ of $Hom(K_{0}({\cal G})^{\otimes n},K_{1}({\cal H}))$.
Now the direct checking shows that the maps 
$$Aut(F^{\otimes n})\longrightarrow Hom(K_{0}({\cal G})^{\otimes n},K_{1}({\cal H}))$$
commutes with the coface operators:
$$\overline{{\partial}^{i}_{n+1}(\alpha)}([X_{1}],...,[X_{n+1}]) =  \left[{\partial}^{i}_{n+1}(\alpha)_{X_{1},...,X_{n+1}}\right] = $$
$$\left\{
\begin{array}{cccc}
i=0, & \left[I_{X_{1}}\otimes\alpha_{X_{2},...,X_{n+1}}\right] &=& [X_{1}]\overline{\alpha}([X_{2}],...,[X_{n+1}])\\
1\leq i\leq n, & \left[\alpha_{X_{1},...,X_{i}\otimes X_{i+1},...,X_{n+1}}\right]
&=& \overline{\alpha}([X_{1}],...,[X_{i}][X_{i+1}],...,[X_{n+1}])  \\ i=n+1, & \left[\alpha_{X_{1},...,X_{n}}\otimes I_{X_{n+1}}\right] &=& 
\overline{\alpha}([X_{1}],...,[X_{n}])[X_{n+1}]
\end{array}
\right\}$$
$$ = {\partial}^{i}_{n+1}(\overline{\alpha})([X_{1}],...,[X_{n+1}])
$$
Thus we have a map 
$${\it h}^{i}(F)\longrightarrow HH^{i}(K_{0}({\cal G}),K_{1}({\cal H}))$$
which asserts to the automorphism $\alpha$ the class $[\overline{\alpha}]$. 

\section{bimodules}
A {\em bimodule} over an associative algebra $R$ is a $k$-module $M$ with the left and right $R$-module structures
$$R\otimes M\rightarrow M \qquad r\otimes m\mapsto rm$$
$$M\otimes R\rightarrow M \qquad m\otimes r\mapsto mr$$
(left and right $k$-module structures coinsides) such that
$$r(ms) = (rm)s\qquad \mbox{for any}\quad r,s\in R, m\in M.$$
A homomorphism ($Hom_{R-R}(M,N)$) from a bimodule $M$ to a bimodule $N$ is a $k$-linear map which preserves both left and right $R$-module structures.
\newline
For example, the endomorphism algebra $End_{R-R}(R\otimes_{k} R)$ of the bimodule $R\otimes_{k} R$ is isomorphic to the algebra $R^{op}\otimes R$ ($R^{op}$ is the algebra with opposite multiplication). Any element $x\in R^{op}\otimes R$ defines the homom
orphism $f_{x}$ (for decomposible $x = r\otimes s$ $f_{x}$ sends $p\otimes q$ from $R\otimes_{k} R$ to $pr\otimes sq$). Conversely the value of an endomorphism $f$ of the bimodule $R\otimes_{k} R$ on the element $1\otimes 1$ lies in $R\otimes_{k} R$.
\newline 
The category $R-Mod-R$ of bimodules over $R$ is an monoidal category with respect to tensor product $\otimes_R$ of bimodules and trivial associativity constraint.
\begin{prop}
The complex of endomorphisms 
$$E_{*}(R-Mod-R) = End(\I_{R-Mod-R}^{\otimes *})$$ 
of the identity functor of the category of bimodules $R-Mod-R$ is isomorphic to the complex 
$Z(R^{\otimes_{k}*+1})$ whose coface maps 
$${\partial}^{i}_{n}:Z(R^{\otimes_{k}n})\longrightarrow Z(R^{\otimes_{k}n+1})$$
are induced by the homomorphisms of algebras 
$${\partial}^{i}_{n}:R^{\otimes_{k}n}\longrightarrow R^{\otimes_{k}n+1}$$
$${\partial}^{i}_{n}(r_{1}\otimes ...\otimes r_{n}) = r_{1}\otimes ...\otimes r_{i}\otimes 1\otimes r_{i+1}\otimes  ...\otimes r_{n}.$$
The isomorphis is realized by two mutially inverse maps:
$$End(\I_{R-Mod-R}^{\otimes n})\longrightarrow Z(R^{\otimes n+1})$$ 
which sends the endomorphism to its specialization on the objects $R\otimes_{k}R,...,R\otimes_{k}R$, and
$$Z(R^{\otimes n+1})\longrightarrow End(\I_{R-Mod-R}^{\otimes n})$$
which associates to an element $r = \sum_{i} r_{0}^{i}\otimes ...\otimes r_{n}^{i}$ the endomorphism $\alpha (r)$ whose specialization on the objects $M_{1},...,M_{n} \in R-Mod-R$ is
$${\alpha (r)}_{M_{1},...,M_{n}}(m_{1}\otimes ...\otimes m_{n}) = \sum_{i} r_{0}^{i}m_{1}r_{1}^{i}\otimes ...\otimes m_{n}r_{n}^{i}.$$
\end{prop}
Proof:
\newline
The bimodule $R\otimes_{k} R$ is a projective generator in the category $R-Mod-R$. Indeed the functor
$$k-Mod\longrightarrow R-Mod-R \qquad V\mapsto R\otimes_{k} V\otimes_{k} R$$ 
is a right adjoint of the forgetfull functor $R-Mod-R\longrightarrow k-Mod$. The adjunction is given by the morphisms 
$$R\otimes_{k} M\otimes_{k} R\rightarrow M \qquad r\otimes m\otimes s\mapsto rms,$$
$$V\rightarrow R\otimes_{k} V\otimes_{k} R \qquad v\mapsto 1\otimes v\otimes 1.$$
Thus the endomorphism algebra $End(\I_{R-Mod-R}^{\otimes n})$ is isomorphic to the centralizer of the subalgebra $End_{R-R}(R\otimes_{k} R)^{\otimes_{R} n}$ in the algebra 
$End_{R-R}(R^{\otimes_{k} n+1})$. 
\newline
It is easy to see that the endomorphism $f$ from $End_{R-R}(R^{\otimes_{k} n+1})$ which commutes with the subalgebra $End_{R-R}(R\otimes_{k} R)^{\otimes_{R} n}$ is defined by its value $f(1\otimes ...\otimes 1) \in R^{\otimes n+1}$. Indeed, since $f$ com
mutes with $f_{r\otimes s}$, we have
$$f(r_{1}\otimes ...\otimes r_{n+1}) = f((r_{1}\otimes 1)\otimes_{R} ...\otimes_{R} (r_{n}\otimes r_{n+1})) =$$ $$f((f_{r_{1}\otimes 1}\otimes_{R}...\otimes_{R} f_{r_{n}\otimes r_{n+1}})((1\otimes 1)\otimes_{R} ...\otimes_{R} (1\otimes 1)) = $$
$$(f_{r_{1}\otimes 1}\otimes_{R}...\otimes_{R} f_{r_{n}\otimes r_{n+1}})(f(1\otimes ...\otimes 1)).$$
Finally, R-polylinearity implice that $f(1\otimes ...\otimes 1)$ lies in $Z(R^{\otimes n+1})$.
\newline
\newline
Since the algebras $End(\I_{R-Mod-R}^{\otimes n})$ are commutative, the cohomology $h^{n}(R-Mod-R) = h^{n}(\I_{R-Mod-R})$ is defined for any $n$ and coincides with the cohomology of the cochain complex associated to the cosimplicial complex of groups of 
invertible elements of the algebras $Z(R^{\otimes n+1})$. In the case of commutative algebra $R$ this complex is called the {\em Amitsur complex} of the algebra $R$ and its cohomology ($HA^{*}(R/k)$) the {\em Amitsur cohomology} \cite{rz}. In particular 
$h^{1}(R-Mod-R)$ is isomorphic to the relative Picard group $Pic(R/k)$ and $h^{2}(R-Mod-R)$ to the relative Brauer group $Br(R/k)$. 
\newline
Using the restriction $h^{n}(R-Mod-R)\rightarrow h^{n}(R-mod-R)$ and the homomorphism from the cohomology of the identity functor to the Hochschild cohomology of the K-theory we can define the following maps
$$HA^{n}(R/k)\rightarrow HH^{n}(K_{0}(R-mod-R),K_{1}(R-mod-R)),$$
for the commutative algebra $R$.
\newline
In the case of a Galois extension of fields $K/k$ with the Galois group $G$ \cite{chr} the previous homomorphism is identical:
$$H^{n}(G,R^{*}) = HA^{n}(R/k)\rightarrow HH^{n}(K_{0}(R-mod-R),K_{1}(R-mod-R)) = H^{n}(G,R^{*}).$$
Indeed, the category of bimodules $R-mod-R$ over Galois extension $K/k$ is semisimple, its simple objects being parametrized by the elements of $G$, whence 
$$K_{*}(R-mod-R) = K_{*}(K)\otimes_{\Bbb Z}{\Bbb Z}[G].$$

\section{modules and comodules over bialgebra}
A {\em bialgebra} \cite{swe} is an algebra $H$ together with a homomoprhisms of algebras 
$$\Delta :H\rightarrow H\otimes H\qquad\mbox{(coproduct)}$$
$$\varepsilon :H\rightarrow k\qquad\mbox{(counit)}$$
for which 
$$(\Delta\otimes I)\Delta = (I\otimes\Delta)\Delta\qquad\mbox{(coassociativity)},$$
$$(\varepsilon\otimes I)\Delta = (I\otimes\varepsilon)\Delta = I\qquad\mbox{(axiom of counit)}.$$
For example, the group algebra $k[G]$ of a group $G$ is a bialgebra, where  
$$\Delta :k[G]\rightarrow k[G]\otimes k[G]\qquad \Delta (g) = g\otimes g,$$
$$\varepsilon :k[G]\rightarrow k\qquad \varepsilon (g) = 1\qquad g\in G.$$
Another example provided by the universal enveloping algebra $U[{\frak g}]$ of the Lie algebra ${\frak g}$
$$\Delta :U[{\frak g}]\rightarrow U[{\frak g}]\otimes U[{\frak g}]\qquad \Delta (l) = l\otimes 1 + 1\otimes l,$$
$$\varepsilon :U[{\frak g}]\rightarrow k\qquad \varepsilon (l) = 0\qquad l\in{\frak g}.$$
\newline 
The coproduct allows to define the structure of $H$-module on the tensor product $M\otimes_{k}N$ of two $H$-modules
$$h*(m\otimes n) = \Delta(h)(m\otimes n)\qquad h\in H, m\in M, n\in N.$$
The coassociativity of coproduct implices that the standart associativity constraint for underlying $k$-modules
$$\varphi :L\otimes (M\otimes N)\rightarrow (L\otimes M)\otimes N \qquad \varphi (l\otimes (m\otimes n) = (l\otimes m)\otimes n$$
induces an associativity for this tensor product.
\newline
The counit defines the structure of $H$-module on the ground ring $k$ which (by counit axiom) is an unit object respectively to the tensor product.
\newline
By another words the category of (left) $H$-modules $H-Mod$ is a monoidal category. It is follows from the definition of tensor product in $H-Mod$ that the forgetful functor $F: H-Mod\rightarrow k-Mod$ is a monoidal functor with trivial monoidal structur
e. 
\begin{prop}
1.The complex of endomorphisms $E_{*}(F) = End(F^{\otimes *})$ of forgetful functor of the category of modules $H-Mod$ over bialgebra $H$ is isomorphic to the bar complex 
$H^{\otimes_{k}*}$ of $H$ (the coface maps 
$${\partial}^{i}_{n}:H^{\otimes_{k}n-1}\longrightarrow H^{\otimes_{k}n}$$
has the form 
$${\partial}^{i}_{n}(h_{1}\otimes ...\otimes h_{n}) = 
\left\{
\begin{array}{cl}
1\otimes h_{1}\otimes ...\otimes h_{n},& i=0\\
h_{1}\otimes ...\otimes\Delta (h_{i})\otimes ...\otimes h_{n},& 1\leq i\leq n\\
h_{1}\otimes ...\otimes h_{n}\otimes 1,& i=n+1
\end{array}
\right.
$$
the codegeneration maps are
$${\sigma}^{i}_{n}(h_{1}\otimes ...\otimes h_{n+1}) = h_{1}\otimes ...\otimes\varepsilon (h_{i})\otimes ...\otimes h_{n+1}))$$
The isomorphis is realized by two mutially inverse maps:
$$End(F^{\otimes n})\longrightarrow H^{\otimes n}$$ 
which sends the endomorphism to its specialization on the objects $H,...,H$, 
\newline
and
$$H^{\otimes n}\longrightarrow End(F^{\otimes n})$$
which associate to the element $x\in H^{\otimes n}$ the endomorphism of multiplying by $x$.
\newline
2.The complex of endomorphisms $E_{*}(H-Mod) = End(\I_{H-Mod}^{\otimes *})$ of the identity functor of the category $H-Mod$ is isomorphic to the subcomplex of bar complex of $H$ which consists of $H$-invariant elements (the subcomplex of centralizers $C_
{H^{\otimes n}}(\Delta (H))$ of the images of diagonal embeddings).
\end{prop}
Proof:
\newline
We will proof more general statement. 
\newline
Let $f:F\rightarrow H$ be a homomorphism of bialgebras (it means that $f$ is a homomorphism of algebras and $\Delta f = (f\otimes f)\Delta$). The restriction functor $f^{*}:H-Mod\rightarrow F-Mod$ is monoidal and the complex of its endomorphisms coincide
s with the subcomplex $C_{H^{\otimes *}}(\Delta (f(F))$ in the bar complex of the bialgebra $H$. The proposition follows from this fact, because the forgetful functor from the category $H-Mod$ is the restriction functor, which corresponds to the unit emb
edding $k\rightarrow H$, and the identity functor is the restriction functor, which corresponds to the identity map $H\rightarrow H$.
\newline
The forgetful functor from the category $H-Mod$ have a right adjoint
$$k-Mod\longrightarrow H-Mod\quad V\mapsto H\otimes_{k} V$$
which endomorphism algebra coincides with an algebra of endomorphisms of $H$-module $H$ and is isomoprhic to $H^{op}$ (the element $h\in H^{op}$ defines the endomorphism $r_{h}$ of right multiplying by $h$).
\newline
Thus the endomorphism algebra $End((f^{*})^{\otimes n})$ is isomorphic to the centralizer of $End_{H}(H)^{\otimes n}$ in the algebra $End_{F}(H^{\otimes n})$.
\newline
It is follows from the direct checking that an element $f$ of this centralizer is defined by its value $f(1\otimes ...\otimes 1) \in H^{\otimes n}$. Indeed, since $f$ commutes with right multiplications by elements of $H$ we have
$$f(h_{1}\otimes ...\otimes h_{n}) = f((r_{h_{1}}\otimes ...\otimes r_{h_{n}})(1\otimes ...\otimes 1)) = (r_{h_{1}}\otimes ...\otimes r_{h_{n}})f(1\otimes ...\otimes 1)$$ 
The composition of endomorphisms from $C_{End_{F}(H^{\otimes n})}(End_{H}(H)^{\otimes n})$ corresponds to the multiplication of its falues in $H^{\otimes n}$
$$fg(1\otimes ...\otimes 1) = f(g(1\otimes ...\otimes 1)) = 
f(\sum_{i} h_{i,1}\otimes ...\otimes h_{i,n}) = $$
$$f(1\otimes ...\otimes 1)(\sum_{i} h_{i,1}\otimes ...\otimes h_{i,n}) = f(1\otimes ...\otimes 1)g(1\otimes ...\otimes 1).$$
Finally, F-polylinearity of $f$ means that $f(1\otimes ...\otimes 1)$ lies in the cetralizer of the image of diagonal embidding $\Delta (F)$ in $H^{\otimes n}$.
\newline
\newline
The first cohomology $h^{1}(F) = z^{1}(F)$ of the forgetful functor $F:H-Mod\rightarrow k-Mod$ coincides with the group of invertible {\em group-like} elements of the bialgebra $H$ \cite{swe}
$$G(H) = \{ g\in H, \Delta (g)=g\otimes g\}.$$
The elements of $h^{2}(F)$ corresponds to the twistings of the coproduct in the Hopf algebra $H$ \cite{dr}. Namely, the adjunction $$\Delta^{x}(h) = x\Delta (h)x^{-1},\ h\in H$$
of the coproduct $\Delta$ by the element from
$$z^{2}(F) = \{x\in (H^{\otimes 2})^* ,\ (1\otimes x)(I\otimes\Delta )(x) = (x\otimes 1)(\Delta\otimes I)(x)\}$$
also satisfies to the coassociativity axiom.
\newline
If $x$ also satisfies to the condition 
$$(\varepsilon\otimes I)(x) = (I\otimes\varepsilon )(x) = 1$$
(which may be attained by the substitution of $x$ by some equivalent to it), then the triple $(H,\Delta^x,\varepsilon)$ is a bialgebra ({\em twisted form} of bialgebra $H$). 
\newline
It is follows from the definitions that the complex of endomorphisms of the twisted form $F(x)$ of forgetful functor $F$ coincides with the bar complex of the twisted form $(H,\Delta^{x},\varepsilon)$ of bialgebra $H$.
\newline
The set $h^{2}(F)$ accepts also another describtion. Namely, $h^{2}(F)$ coincides with the set of Galois $H$-module coalgebras $Gal_{k-coalg}(H)$, which is isomorphic to $H$ as $H$-modules \cite{ul}.
\newline
The {\em $H$-module coalgebra} is a coalgebra $L$ which coproduct is a homomorphism of $H$-modules. The $H$-module coalgebra $L$ is {\em Galois} \cite{cs} if the map
$$H\otimes L\rightarrow L\otimes L,\qquad h\otimes l\mapsto (h\otimes 1)\Delta (l)$$
is an isomorphism.
\newline
We have a map 
$$h^{2}(F)\rightarrow Gal_{k-coalg}(H),$$
which sends the class of cocycle $x\in z^{2}(F)$ to the calss of coalgebra $(H,\Delta_{x})$, where  $\Delta_{x}(h) = \Delta (h)x^{-1},\ h\in H$.
\newline
The first cohomology $h^{1}(H-Mod) = z^{1}(\I_{H-Mod})$ of the identity functor coincides with the centre of the group of invertible group-like elements of the bialgebra $H$ 
$$Z(G(H)) = \{ g\in Z(H)^{*}, \Delta (g)=g\otimes g\}.$$
\newline
The kernel of the map of punctured sets 
$$h^{2}(H-Mod)\rightarrow h^{2}(F)$$
(the stabilizer of the action of the group $h^{2}(H-Mod)$ on the set $h^{2}(F)$) coincides with the kernel of natural homomorphism
$$Out_{bialg}(H)\rightarrow Out_{alg}(H)$$
from the group $Out_{bialg}(H)$ of outer automorphisms of bialgebra $H$ (the factor of the group of automorphism of bialgebra $H$ modulo conjugations by group-like elements of $H$) to the group $Out_{alg}(H)$ of outer automorphisms of the algebra $H$. In
deed, an element of second kernel is a conjugation by an element $h\in H^*$ such that
$y = \Delta (h)(h\otimes h)^{-1}$ lies in the centralizer $C_{H^{\otimes 2}}(\Delta (H))$. Hence, $y$ represents the element of $h^{2}(H-Mod)$ which trivialises in $h^{2}(F)$.
\newline
For example, in the case of the group bialgebra $k[G]$ the first cohomology $h^{1}(F)$ of forgetful functor is isomorphic to the group $G$. The second cohomology $h^{2}(F)$ coincides with the set $Gal_{k}(G)$ of Galois $G$-extension of $k$. The kernel of
 the map
$$h^{2}(k[G]-Mod)\rightarrow h^{2}(F)$$
coincides (in the case of finite group $G$ and the field $k$ such that $char(k) \not\mid |G|$) with the factorgroup of {\em locally-inner} automorphisms of $G$ (an automorphisms which not moves the conjugacy classes) modulo inner.
\newline
The cohomology of forgetful and identity functors in the case of the universal enveloping algebra $U[{\frak g}]$ of the Lie algebra (over the fild $k$) are trivial. It follows from the absence of non-trivial invertible elements in $U[{\frak g}]$, which c
an be deduced from Poincare-Birkhoff-Witt theorem. But non-trivial invertible elements appears if we extend the scalars (ground ring $k$). For example, the tangent cohomology of forgetful and identity functors are non-trivial in general. 
\newline
It was shown in \cite{dr} that the natural homomorphism
$$\Lambda^{*}{\frak g} = \Lambda^{*}H^{1}(F)\rightarrow H^{*}(F),$$
induced by the multiplication in $H^{*}(F)$ is an isomorphism. The first approximation of the tangent cone in $H^{2}(F)$ is given by the equation
$$[[\alpha ,\alpha ]] = 0,$$
where $[[\ ,\ ]]$ is a component of the Lie superalgebra structure on $\Lambda^{*}{\frak g}$ \cite{dri}
$$[[x_{1}\wedge ...\wedge x_{s},y_{1}\wedge ...\wedge y_{t}]] = $$
$$\sum_{i,j}(-1)^{i+j}[x_{i},y_{j}]\wedge x_{1}\wedge ...\wedge\widehat{x_{i}}\wedge ...\wedge x_{s}\wedge y_{1}\wedge ...\wedge\widehat{y_{j}}\wedge ...\wedge y_{t},$$
where $\widehat{z}$ means that $z$ not accours in the product .
\newline
It follows from the result of \cite{ek} that the tangent cone in $H^{2}(F)$ coincides with its first approximation.
\newline
In particular, for any Cartan subalgebra ${\frak h}$ the space $\Lambda^{2}{\frak h}$ (the so-called infinitesimal {\em Sudbery family}) lies in the first approximation and in the tangent cone \cite{ggs,sad}.
\newline
It can be deduced from the description of $H^{*}(F)$ that the tangent cohomology of the identity functor $H^{*}(U[{\frak g}]-Mod)$ coincides with $\frak g$-invariant elements in 
$\Lambda^{*}{\frak g}$.
\newline
\newline
There is another type of monoidal categories which can be  associated to the bialgebras, namely categories of comodules.
\newline
A (right) {\em comodule} over a bialgebra $H$ \cite{swe} is a $k$-module $M$ with the homomorphism $\psi:M\rightarrow M\otimes H$ such that 
$$(I\otimes\Delta )\psi = (\psi\otimes I)\psi\qquad\mbox{(coassociativity of comodule structure)},$$
$$(I\otimes\varepsilon )\psi = I\qquad\mbox{(counitarity)}.$$
For example, the ground ring $k$ has the comodule structure 
$$k\rightarrow k\otimes H = H,\qquad c\mapsto c\otimes 1=c1.$$
Another example provides by the bialgebra $H$ with the coproduct $\Delta$ regarded as a comodule structure.
\newline
A {\em morphism of comodules} $(M,\psi)$ and $(N,\phi)$ is a $k$-linear map $f:M\rightarrow N$ which preserves the comodule structure: 
$$\phi f = (f\otimes I)\psi.$$
We will denote the set of morphisms of comodules by $Hom^{H}(M,N)$. For example, the endomorphism algebra $End^{H}(H)$ of comodule $H$ is isomorphic to the $Hom_{k}(H,k)$ ({\em dual algebra} to $H$) with the multiplication ({\em convolution})
$$p*q = (p\otimes q)\Delta,\qquad p,q\in Hom_{k}(H,k).$$
Indeed, an endomorphism $f\in End^{H}(H)$ defines a map $p_{f}=\varepsilon f\in Hom_{k}(H,k)$. Conversely, a linear map $p:H\rightarrow k$ defines an endomorphism $f_{p} = (p\otimes I)\Delta$.
\newline
The product $\mu :H\otimes H\rightarrow H$ of the bialgebra $H$ allows to define a comodule structure on the tensor product (over $k$) of any two comodules $(M,\psi)$ and $(N,\phi)$
$$\gamma:M\otimes N\rightarrow M\otimes N\otimes H,\qquad \gamma = (I_{M\otimes N}\otimes\mu)t_{H,N}(\psi\otimes\phi).$$
The associativity of the product $\mu$ implices that the (trivial) assocociativity constraint of underling $k$-modules is an associativity for tensor product $\otimes_{k}$ in the category of comodules $Comod-H$.
It is emmideatly follows from the definitions that the comodule $k$ is a unit object of this category and the forgetful functor $F:Comod-H\rightarrow k-Mod$ is a monoidal functor with trivial monoidal structure. 
\begin{prop}
1.The complex of endomorphisms $E_{*}(F) = End(F^{\otimes *})$ of forgetful functor of the category of comodules $Comod-H$ over bialgebra $H$ is isomorphic to the cobar complex 
$Hom_{k}(H^{\otimes*},k)$ of $H$ (with coface maps 
$${\partial}^{i}_{n}:Hom_{k}(H^{\otimes n-1},k)\longrightarrow Hom_{k}(H^{\otimes n},k)$$
which are induced by the homomorphisms
$d^{i}_{n}:H^{\otimes n}\longrightarrow H^{\otimes n-1}$ 
$$d^{i}_{n}(h_{1}\otimes ...\otimes h_{n}) = 
\left\{
\begin{array}{cl}
\varepsilon (h_{1})\otimes ...\otimes h_{n},& i=0\\
h_{1}\otimes ...\otimes h_{i}h_{i+1}\otimes ...\otimes h_{n},& 1\leq i\leq n\\
h_{1}\otimes ...\otimes\varepsilon (h_{n}),& i=n+1
\end{array}
\right.
$$
codegeneration maps are induced by the homomorphisms $s^{i}_{n}:H^{\otimes n}\rightarrow H^{\otimes n+1}$ 
$$s^{i}_{n}(h_{1}\otimes ...\otimes h_{n}) = h_{1}\otimes ...\otimes\Delta (h_{i})\otimes ...\otimes h_{n}))$$
The isomorphism is realized by two mutially inverse maps:
$$End(F^{\otimes n})\longrightarrow Hom(H^{\otimes n},k)$$ 
which sends the endomorphism to its specialization on the objects $H,...,H$, and
$$Hom(H^{\otimes n},k)\longrightarrow End(F^{\otimes n})$$
which associate to the element $x\in Hom(H^{\otimes n},k)$ the endomorphism of convolution with $x$.
\newline
2.The complex of endomorphisms $E_{*}(Comod-H) = End(\I_{Comod-H}^{\otimes *})$ of the identity functor of the category $Comod-H$ is isomorphic to the subcomplex of the cobar complex of $H$ which consists of $H$-coinvariant elements (the subcomplex of ce
ntralizers $C_{Hom(H^{\otimes n},k)}(\delta (Hom(H,k)))$ of the images of diagonal embeddings of $Hom(H,k)$).
\end{prop}
Proof:
\newline
As in that case of modules over bialgebra we can proof more general statement. Let $f:H\rightarrow F$ is a homomorphism of bialgebras. The corestriction functor 
$$f_{*}:H-Mod\rightarrow F-Mod,\qquad f_{*}(M,\psi) = (M,(I\otimes f)\psi)$$ 
is a monoidal and the complex of its endomorphisms coincides with the subcomplex $C_{Hom(H^{\otimes *},k)}(\delta (f(Hom(F,k)))$ in the cobar complex of the bialgebra $H$. The proposition follows from this fact, because the forgetful functor from the cat
egory $H-Mod$ is the corestriction functor which corresponds to the counit $H\rightarrow k$ and the identity functor is the corestriction functor which corresponds to the identity map $H\rightarrow H$.
\newline
The forgetful functor from the category $Comod-H$ have a left adjoint
$$k-Mod\longrightarrow Comod-H\quad V\mapsto V\otimes_{k} H$$
which endomorphism algebra coincides with the algebra of endomorphisms of comodule $H$.
\newline
Thus the endomorphism algebra $End((f^{*})^{\otimes n})$ is isomorphic to the centralizer of $End^{H}(H)^{\otimes n}$ in the algebra $End^{F}(H^{\otimes n})$.
\newline
It is follows from the direct checking that the element $\alpha$ of this centralizer is defined by the map $p_{\alpha} = (\varepsilon\otimes ...\otimes\varepsilon)f \in Hom(H^{\otimes n})$.
It is easy to see that the composition of the endomorphisms $\alpha$ and $\beta$ from $C_{End_{F}(H^{\otimes n})}(End_{H}(H)^{\otimes n})$ corresponds to the convolution of $p_{\alpha}$ and $p_{\beta}$. 
Finally, F-polylinearity of $\alpha$ means that $p_{\alpha}$ lies in the cetralizer of the image of diagonal embidding $\delta (Hom(F,k))$ in $Hom(H^{\otimes n},k)$.
\newline
\newline
For a {\em cocommutative} bialgebra $H$ (cocommutativity means that $t\Delta = \Delta$, where $t$ is a permutation of factors) the cobar complex is $\infty$-comutative. In particular, the cohomology $h^{n}(F)$ of forgetful fuctor $F:Comod-H\rightarrow K-
Mod$ coincides with the cohomology of identity functor $h^{n}(Comod-H)$ and is well defined for any $n$.
\newline
This cohomology was considered by Sweedler \cite{sw:c}. He proved that in the case of group bialgebra of the group $G$ the cohomology of forgetfull functor is isomorphic to the group cohomology of $G$ with coefficients in the (trivial $G$-module) of inve
rtible elements of the groun ring $k$
$$h^{n}(F_{Comod-k[G]}) \simeq H^{n}(G,k^{*})$$
and in the case of universal enveloping algebra $U[{\frak g}]$ of the Lie algebra $\frak g$ the cohomology of forgetfull functor coincides with the cohomology of Lie bialgebra $\frak g$ with coefficients in the groun ring $k$
$$h^{n}(F_{Comod-U[{\frak g}]}) \simeq H^{n}({\frak g},k).$$
\newline
Let us note that for the category of comodules over group algebra of the group $G$ the homomorphism from the cohomology of the forgetful functor to the Hochschild cohomology of K-theory
$$H^{n}(G,k^{*}) = h^{n}(F_{Comod-k[G]})\rightarrow HH^{n}(K_{0}(Comod-k[G]),K_{0}(k-Mod)) = H^{n}(G,k^{*})$$
is identical. 
\newline
It follows from the fact that the category $Comod-k[G]$ (the category of $G$-graded vector spaces) is 
a semisimple and all its simple objects are invertible (and parametrized by $G$).

\section{free tensor category. unitary R-matrices}
Let now $k$ is a field. 
The free $k$-linear (abelian) tensor category ${\cal T}_k$ \cite{dm} is a cartesian product
\newcommand{\optimes}{\mathop{\times}\limits}
$$\optimes_{n\geq 0}k[S_{n}]-mod$$
of the categories of finite-dimensional representations of the symmetric groups $S_{n}$. 
\newline
For the field of characteristic zero ${\cal T}_k$ is a semisimple category whose simple objects corresponds to the partitions. We will denote by $S^{\mu}X$ the simple object of ${\cal T}_k$ labeled by the partition $\mu$. Specifically, the object $1 = S^
{(0)}X$ corresponding to the (trivial) partition of $0$ is a unit object in ${\cal T}_k$ and the object $X = S^{(1)}X$ which corresponds to the (trivial) partition of $1$ is a tensor generator in ${\cal T}_k$ (any object of ${\cal T}_k$ is a direct summu
nd of some sum of tensor powers of $X$). Note that 
$$Hom_{\cal T}(X^{\otimes n},X^{\otimes m}) = 
\left\{
\begin{array}{cl}
k[S_{n}],& n=m\\
0,& n\not=m
\end{array}
\right.$$
It was pointed by Yu.I.Manin \cite{man} that a monoidal functor from the category ${\cal T}_k$ to the category of $k$-vector spases $k-Mod$ is nothing as a unitary solution of the {\em quantum Yang-Baxter} equation ({\em unitary quantum $R$-matrix}), i.e
. an automorphism $R\in Aut(V^{\otimes 2})$ of tensor square of some vector space $V$
such that 
$$R_{1}R_{2}R_{1} = R_{2}R_{1}R_{2},\qquad R^2 = 1,$$ 
where $R_{1} = R\otimes I_V$ and $R_{2} = I_V\otimes R$ are an automorphisms of tensor qube of $V$.
\newline
Indeed, the value $F(\tau)$ of the monoidal functor $F$ on the automorphism $\tau\in Aut(X^{\otimes 2})$ (the generator of $S_{2}$) is a unitary quantum $R$-matrix on the vector space $V = F(X)$. The condition $R^2 = 1$ is obvious and quantum Yang-Baxter
 equation follows from Coxeter relation in $S_{3}$. The pair $(F(X),F(\tau))$ is defined the functor $F$ because the category ${\cal T}_k$ is generated by the object $X$ and the morphism $\tau$.
\newline
It is natural to ask when the functors $F_R$ and $F_S$, corresponding to the quantum $R$-matrices $R$, and $S$ are isomorphic as functors (are twisted forms of each other as monoidal functors). 
\newline
It is easy to see that two additive functors $F$ and $G$ between semisimple categories ${\cal A}$ and ${\cal B}$ are isomorphic if and only if they induces the same homomorphism between Grothendieck groups 
$$K_{0}(F)=K_{0}(G):K_{0}({\cal A})\rightarrow K_{0}({\cal B}).$$ 
In the case of characteristic zero the Grothendieck ring of the free tensor category ${\cal T}_k$ coincides with the free $\lambda$-ring (generating by one element $x=[X]$) and is isomorphic (as a ring) to the ring of polinomials of infinitely many varia
bles \cite{knu}
$$K_{0}({\cal T}_k) = {\Bbb Z}[x,\lambda^{2}x,\lambda^{3}x,...],$$
where $\lambda^{n}x = [\Lambda^{n}X]$ is a class of the  simple object $\Lambda^{n}X$ corresponding to the non-trivial one-dimensional representation of $S_{n}$. The homomorphism $f$ from $K_{0}({\cal T}_k)$ to ${\Bbb Z}=K_{0}(k-mod)$ is defined by its v
alues $f(\lambda^{n}x)$ or by the formal series $H_{f}(t) = \sum_{n\geq 0}f(\lambda^{n}x)t^n$ ({\em Hilbert series} of $f$). This term comes from the fact that $\Lambda^{*}X = \bigoplus_{n\geq 0}\Lambda^{n}X$ is a graded  algebra in (appropriate extensio
n of) the category ${\cal T}_k$ ({\em exterior algebra of $X$}). For any monoidal functor $F:{\cal T}_k\rightarrow k-mod$ $F(\Lambda^{*}X)$ is a graded algebra and $H_{K_{0}(F)}(t)$ coincides with its Hilbert series. 
\begin{prop}
Two $R$-matrices $R$ and $S$ (over the field $k$ of characteristic zero) defines isomorphic monoidal functors $$F_{R},F_{S}:{\cal T}_k\rightarrow k-mod$$
if and only if its Hilbert series $H_{R}(t)=H_{F_{R}}(t), H_{S}(t)=H_{F_{S}}(t)$ coincides.
\end{prop}
The functor $F:{\cal G}\rightarrow {\cal H}$ between monoidal categories is called {\em quasimonoidal} if it preservs the tensor product 
$$F(X\otimes Y)\simeq F(X)\otimes F(Y),\qquad\mbox{for any}\quad X,Y\in{\cal G}.$$
For the quasimonoidal functor $F:{\cal T}_k \rightarrow k-mod$ the nonegativity of values of $K_{0}(F)$ on the classes of simple objects $S^{\mu}X$ implices some conditions on the Hilbert series $H_{K_{0}(F)}(t)$.
\begin{prop}
Let $k$ is an algebraically closed field of characteristic zero. The homomorphism $f:K_{0}({\cal T}_k)\rightarrow {\Bbb Z}$ is a class $K_{0}(F)$ of some quasimonoidal functor $F:{\cal T}_k\rightarrow k-mod$ if and only if there are real positive
$a_{i}, b_{j},\ i=1,...,n,\ j=1,...,m$ such that
$$H_{f}(t) = \frac{\prod_{i=0}^{n}(1+a_{i}t)}{\prod_{i=0}^{m}(1-b_{i}t)}.$$
\end{prop}
Proof:
\newline
It is not hard to verify that the homomorphism $F:K_{0}({\cal A})\rightarrow K_{0}({\cal B})$ between Grothendieck rings of semisimple $k$-linear monoidal categories with finite-dimensional spaces of morhisms is a class of some quasimonoidal functor $F:{
\cal A}\rightarrow {\cal B}$ iff the values of $f$ on the classes of simple objects of ${\cal A}$ are non-negative (has non-negative coefficients in expressions via the classes of simple objects of ${\cal B}$).
\newline
The class $s^{\mu}x$ of the object $S^{\mu}X$ as an element of $K_{0}({\cal T}_k)$ can be written as polinomial of $\lambda^{n}x$. Namely \cite{knu,mcd} 
$$s^{\mu}x = det\left(\lambda^{{\mu}_{i}'-i+j}x\right)_{1\leq i,j\leq n},$$
where ${\mu}' = ({\mu}_{1}',...,{\mu}_{m}')$ is a dual partition for partition $\mu$ (${\mu}_{i}'=|\{j,{\mu}_{j}\geq i\}|$). Hence the non-negativity of values of the homomorphism $f:K_{0}({\cal T}_k)\rightarrow {\Bbb Z}$ on $s^{\mu}x$ for all partitions
 $\mu$ means total positivity of the sequence $\{f(\lambda^{n}x)\}$ (by definition {\em total positivity} \cite{kar} of the sequence $\{a_{n}\}$ is  non-negativity of all finite minors of the infinite Tepliz matrix
$$\left(\begin{array}{cccc}
a_{0} & a_{1} & a_{2} & \cdots \\
0     & a_{0} & a_{1} & \cdots \\
0     & 0     & a_{0} & \cdots \\
\vdots&\vdots &\vdots & \ddots 
\end{array}\right)$$
It is known \cite{kar} that the generating function $\sum_{n\geq 0}a_{n}$ of the total positive sequence $\{a_{n}\}$ have the form 
$$e^{ct}\frac{\prod_{i=0}^{\infty}(1+a_{i}t)}{\prod_{i=0}^{\infty}(1-b_{i}t)},$$
for some real positive $a_{i}, b_{j}, c$ such that $\sum_{i}a_{i}+b_{j} < \infty$. 
\newline
Hence the generating function is meromorphic in $\{t, |t|\leq 1\}$ and has finitly many zeros and poles in this circle. One theorem of Salem \cite{sal} implices that this function is rational.
\newline
\newline
The Hilbert series of $R$-matrix can be expressed in terms of some its symbolic invariants. 
\begin{prop}\cite{man:sem}
Consider the series 
$$\Psi_{R}(t) = (dimV)t + \sum_{n\geq 2}tr(R_{1}...R_{n-1})t^n$$ 
for the $R$-matrix $R$. Then
$$H_{R}(t) = exp\left(\int \frac{\Psi_{R}(-t)}{t}dt \right)$$
\end{prop}
Proof:
\newline
The object $\Lambda^{n}X$ is a image of the idempotent 
$$P_{\Lambda^{n}X} = \frac{1}{n!}\sum_{\sigma\in S_{n}}sign(\sigma )\sigma.$$
Hence the dimension of $F_{R}(\Lambda^{n}X)$ can be expressed in terms of the traces of some products of R-matrix $R$
$$dimF_{R}(\Lambda^{n}X) = Tr(F_{R}(P_{\Lambda^{n}X})) = 
\frac{1}{n!}\sum_{\sigma\in S_{n}}sign(\sigma )Tr(R_{\sigma}),$$ 
where $R_{\sigma} = F_{R}(\sigma)$. Since any permutation $\sigma$ is similar to the product of independent cycles $\tau_{1}...\tau_{k}$ (of lenght $l_{1},...,l_{k}$ respectively), then 
$$Tr(R_{\sigma}) = \psi^{l_{1}}...\psi^{l_{k}},\quad\mbox{where}\ \psi^{l}_{R} = tr(R_{1}...R_{l-1})\ (dimV, l=1).$$
The number of permutations in $S_{n}$ which has the cyclic structure corresponding to the partition $\mu = (\mu_{1},...,\mu_{m}) \in {\cal P}_{n}$ of the number $n$ ($\sum i\mu_i = n$) equals 
$\frac{n!}{\prod_i \mu_{i}!i^{\mu_i}}$. The sign of these partitions is $(-1)^{n+m}$, where $m = \sum\mu_i$. Hence
$$Tr(F_{R}(P_{\Lambda^{n}X})) = \sum_{\mu\in{\cal P}_{n}} (-1)^n \prod_i \left( \frac{-\psi^{l}_{R}}{i} \right)^{\mu_i}\frac{1}{\mu_{i}!}$$
and 
$$H_{R}(t) = \sum_{n\geq 0} (-t)^{n}\sum_{\mu\in{\cal P}_{n}} \prod_i \left(\frac{-\psi^{l}_{R}}{i}\right)^{\mu_i}\frac{1}{\mu_{i}!}.$$
On the other side,
$$exp\left(\int \frac{\Psi_{R}(-t)}{t}dt\right) = \sum_{n\geq 0} \frac{1}{n!}\left(\sum_{i\geq 0} (-t)^i\frac{\psi^{i}_{R}}{i}\right)^n =$$
$$\sum_{\mu\in{\cal P}_{n}} (-1)^n \prod_i \left(\frac{-\psi^{l}_{R}}{i}\right)^{\mu_i}\frac{1}{\mu_{i}!}.$$
\newline
\newline
\begin{exam}{\rm
If $R$ is a unitary $R$-matrix, then $S=-R$ is also unitary $R$-matrix. Since 
$$tr(S_{1}...S_{n-1}) = (-1)^{n-1}tr(R_{1}...R_{n-1}),$$
then
$$\Psi_{S}(t) = dimVt - \sum_{n\geq 2}tr(R_{1}...R_{n-1})(-t)^n = -\Psi_{R}(-t)$$ 
and
$$H_{S}(t) = exp\left(\int \frac{\Psi_{S}(-t)}{t}dt \right) = exp\left(-\int \frac{\Psi_{R}(t)}{t}dt \right) = H_{R}(-t)^{-1}.$$
So the set of $R$-matrices with Hilbert series $H(t)$ is isomorphic to the set of $R$-matrices with Hilbert series $H(-t)^{-1}$.
}\end{exam}
\begin{exam}{\rm
Let $A$ be a commutative algebra,  $a\in A\otimes A$  an element such that $at(a)=1$ ($t:A\otimes A\rightarrow A\otimes A$ is a permutation of the factors) and $M$  a $A$-module. Then 
$$R = R(A,a,M) = L(a)t \in Aut(M\otimes M)$$
is a unitary $R$-matrix on the vector space $M$ (here $L(a)$ is a operator of left multiplying by $a$). Indeed,
$$R_{1}R_{2}R_{1} = L(a_{12})t_{1}L(a_{23})t_{2}L(a_{12})t_{1} =$$
$$L(a_{12}t_{1}(a_{23})t_{1}t_{2}(a_{12}))t_{1}t_{2}t_{1} = 
L(a_{12}a_{13}a_{23})t_{1}t_{2}t_{1}.$$
On the other hand
$$R_{2}R_{1}R_{2} = L(a_{23})t_{2}L(a_{12})t_{1}L(a_{23})t_{2} =$$
$$L(a_{23}t_{2}(a_{12})t_{2}t_{1}(a_{23}))t_{2}t_{1}t_{2} = 
L(a_{23}a_{13}a_{12})t_{2}t_{1}t_{2}.$$
Hence the quantum Yang-Baxter equation for $R$ follows from commutativity of $A^{\otimes 3}$ and Coxeter relation for the permutation $t$. The unitarity of $R$ follows from the condition $at(a)=1$.
\newline
Define some invariant of elements $a\in A\otimes A$ for which  $at(a)=1$. Since algebra $A$ is commutative, then the multiplication $\mu:A\otimes A\rightarrow A$ is a homomorphism of algebras and
$$1 = \mu (at(a)) = \mu (a)\mu (t(a)) = (\mu (a))^2.$$   
So $\mu (a) = \pm 1$.
\newline
Now we can calculate the Hilbert series of $R$. Since
$$R_{1}...R_{n-1} = L(a_{12})t_{1}...L(a_{n-1n})t_{n-1} =$$ 
$$L(a_{12}t_{1}(a_{23})...t_{1}...t_{n-1}(a_{n-1n}))t_{1}...t_{n-1} = L(a_{12}...a_{1n})t_{1}...t_{n-1},$$
then
$$tr(R_{1}...R_{n-1}) = tr(L(a_{12}...a_{1n})t_{1}...t_{n-1}) =$$
$$tr(\mu(a_{12}...a_{1n})t_{1}...t_{n-1}) = (\mu (a))^{n-1}dimM.$$
Hence 
$$\Psi_{R}(t) = dimM\sum_{n\geq 1}(\mu (a))^{n-1}(-t)^n = \frac{dimMt}{1-\mu (a)t}$$
and
$$H_{R}(t) = exp\left(\int \frac{dimM}{1-\mu (a)t}dt\right) = (1+\mu (a)t)^{\mu (a)dimM}.$$
}\end{exam}

\section{free tensor category. Hecke algebras}
This section is devoted to the investigation of the monoidal structures on the category ${\cal T}_k$.
\newline
We begin with the classification of (quasi)monoidal autoequivalences of the category ${\cal T}_k$.
\begin{prop}\label{hom}
Let $f$ is a graded homomorphism $K_{0}({\cal T}_{k})\rightarrow K_{0}({\cal T}_{k})$ for which
$$f(x)=x,\ f(\lambda^{2}x)=\lambda^{2}x,\ \mbox{and}$$
$$f(s^{\mu}x)=\sum n_{\mu,\nu}s^{\nu}x,\ \mbox{where}\ n_{\mu,\nu}\geq 0.$$
Then $f$ is identical. 
\end{prop}
Proof:
\newline
Let us prove the lemma by the induction on the degree.
\newline
Suppose that $f(s^{\mu}x)=s^{\mu}x$ for any partition  $\mu\in{\cal P}_m$ where $m\leq n$.
\newline
To prove that $f(s^{\mu}x)=s^{\mu}x$ for $\mu\in{\cal P}_{n+1}$ it is sufficient to show that $f(\lambda^{n+1}x)=\lambda^{n+1}x$.
It follows from the partial case of the Littlewood-Richardson formula
$$x\cdot s^{\mu}x = \sum_{\mu\subset\nu,|\nu\setminus\mu|=1} s^{\nu}x$$
that
$$f(\lambda^{n+1}x) + f(s^{(n,1)}x) = f(x\cdot\lambda^{n}x) = x\lambda^{n}x = \lambda^{n+1}x + s^{(n,1)}x.$$
The summand $f(s^{(n,1)}x)$ can not equal to zero, because in this case 
$$s^{2}x\cdot\lambda^{n}x = f(s^{2}x)f(\lambda^{n}x) = f(s^{2}x\cdot\lambda^{n}x) = f(\sum_{\nu\subset (n,1)}s^{\nu}x) = 0.$$ 
By another side the expression of $f(s^{(n,1)}x)$ can not contain $\lambda^{n+1}x$, because the right (and hence left) side of the equality
$$f(s^{(n,1)}x) + f(s^{(n-1,1,1)}x) + f(s^{(n-1,2)}x) = f(x\cdot s^{(n-1,1)}x) = x\cdot s^{(n-1,1)}x = $$
$$s^{(n,1)}x + s^{(n-1,1,1)}x + s^{(n-1,2)}x$$
not contain it.
\newline
Hence $f(s^{(n,1)}x) = s^{(n,1)}x$ and $f(\lambda^{n+1}x)=\lambda^{n+1}x$.
\newline
\newline
Using the previous proposition we can describe the (quasi)monoidal autoequivalences of ${\cal T}_{k}$.
\begin{prop}
There is unique non-identical (quasi)monoidal autoequivalence $F$ of the category ${\cal T}_{k}$ which is defined by the seting
$$F(X) = X,\qquad F(\tau ) = -\tau.$$
\end{prop}
Proof:
\newline
Let $G$ is a quasimonoidal autoequivalence of the category ${\cal T}_{k}$.
\newline
Firstly let us show that the automorphism $K_{0}(G)$ of the Grothendieck ring $K_{0}({\cal T}_{k})$ either identical or coincides with $K_{0}(F)$ which sends $s^{\mu}x$ to $s^{{\mu}'}x$ (where ${\mu}'$ is a dual partition to ${\mu}$).
\newline
It is easy to see that the quasimonoidal autoequivalence must send the simple object to the simple. In particular, $G(X)=X$, since $X$ is unique object of ${\cal T}_{k}$ whose square is a sum of two simple objects. For $G(\Lambda^{2}X)$ there are two pos
ibilites $\Lambda^{2}X$ and $S^{2}X$. If $G(\Lambda^{2}X) = S^{2}X$ we may replace $G$ by $FG$, so we can assume that $G(\Lambda^{2}X) = \Lambda^{2}X$. Hence $K_{0}(G)$ is identity by the previous proposition.
\newline
The case of monoidal autoequivalence is more easy. From the previous consideration follows that such autoequivalens sends the generator $X$ to itself and the direct checking shows that there are only two solutions in $End(X^{\otimes 2})$ of the equations

$$t^2 = 1,\qquad t_{1}t_{2}t_{1} = t_{2}t_{1}t_{2}.$$
Namely, the standart $t=\tau$ and $t=-\tau$.
\newline
\newline
Now we will construct of the map
$$a:Ms({\cal T}_{k})\rightarrow k.$$
Let $\psi$ some associativity constraint of the category ${\cal T}_{k}$.
Consider two homomorphisms of algebras
$$f_{1},f_{2}:End_{\cal T}(X^{\otimes 2})\rightarrow End_{\cal T}(X^{\otimes 3}),\qquad f_{1}(g)=g\otimes I,\ f_{2}(g)=\psi^{-1}(I\otimes g)\psi ,$$
here $X^{\otimes 3}$ means $X\otimes (X\otimes X)$. 
\newline
As an algebras
$$End_{\cal T}(X^{\otimes 2}) = End_{\cal T}({\Lambda}^{2}X)\oplus End_{\cal T}(S^{2}X) \simeq k\oplus k$$
and
$$End_{\cal T}(X^{\otimes 3}) = End_{\cal T}({\Lambda}^{3}X) \oplus End_{\cal T}(S^{(2,1)}X \oplus S^{(2,1)}X) \oplus End_{\cal T}(S^{3}X) \simeq k \oplus M_{2}(k) \oplus k,$$
here $(2,1)$ is a partition of the number 3.
\newline
Since 
$${\Lambda}^{\otimes 2}X \otimes X \simeq {\Lambda}^{\otimes 3}X \oplus S^{(2,1)}X$$
the images $p_{i}=f_{i}(p)$ of the projector $p\in End_{\cal T}(X^{\otimes 2})$ on the ${\Lambda}^{\otimes 2}X$ have the decompositions $(1,P_{i},0)\in k \oplus M_{2}(k) \oplus k$.
In addition the projectors $P_{i}\in M_{2}(k)$ has rank 1, hence
$$P_{1}P_{2}P_{1} = aP_{1},\quad P_{2}P_{1}P_{2} = aP_{2}\qquad\mbox{for some} a=a(\psi )\in k.$$
\begin{exam}{\rm
A {\em Hecke algebra} $H_{n}(q),\ q\in k$ is an algebra with generators $t_{i},\ i=1,...,n-1$ and defining relations
$$t_{i}t_{i+1}t_{i}=t_{i+1}t_{i}t_{i+1},\qquad i=1,...,n-2,$$
$$t_{i}t_{j}=t_{j}t_{i}, \qquad |i-j|>2,$$
$$(t_{i}+1)(t_{i}-q)=0,\qquad i=1,...,n-1.$$
It is known \cite{dj,mur} that in the case when $q$ is not a root of unity the algebras $H_{n}(q)$ are semisimple for any $n$ and its irreducible representations are parametrized by the partitions.
\newline
A {\em free heckian category generated by one object} ${\cal T}_{k,q}$ is a cartesian product
$${\cal T}_{k,q} = \times_{n\geq 0} H_{n}(q)-mod,$$ 
(here $H_{0}(q)=H_{1}(q)=k$) with the tensor product, which is induced by the homomorphisms of algebras
$$H_{n}(q)\otimes H_{m}(q) \rightarrow H_{n+m}(q),\qquad t_{i}\otimes 1\mapsto t_{i},\ 1\otimes t_{i}\mapsto t_{n+i}.$$
We will denote by $1,X$ the objects of ${\cal T}_{k,q}$ which corresponds to the unique one-dimensional representations of $H_{0}(q),H_{1}(q)$ and by $\Lambda^{n}_{q}X$ the objects corresponding to the one-dimensional representations of $H_{n}(q)$ which 
sends $t_{i}$ to $-1$.  
\newline
The term free heckian category is explained by the fact that the category ${\cal T}_{k,q}$ is generated (as a monoidal category) by the object $X$ and the automorphism $t\in End_{{\cal T}_{q}}(X^{\otimes 2})$
such that $(t+1)(t-q)=0$ and $t_{1}t_{2}t_{1}=t_{2}t_{1}t_{2}$, where as usually $t_{1}=t \otimes 1,\ t_{2}=1 \otimes t \in End_{{\cal T}_{q}}(X^{\otimes 3})$.
\newline
It is known that (for the case of characteristic zero) the Littlewood-Richardson coefficients of the category ${\cal T}_{k,q}$ coincides with the Littlewood-Richardson coefficients of the category ${\cal T}_{k}$. Hence the monoidal category ${\cal T}_{k,
q}$ is the category ${\cal T}_{k}$ with another associativity constraint ${\varphi}_q$.
\newline
Let us note \cite{dj,mur} that ${\Lambda}^{n}_{q}X$ is a image of the projector 
$$p_{\Lambda^{n}X} = \frac{1}{n!}\sum_{\sigma\in S_n} (-1)^{l(\sigma )}t_{\sigma},$$
(here $l(\sigma )=m$ and $t_{\sigma}=t_{i_1}...t_{i_m}$, if $\sigma=\tau_{i_1}...\tau_{i_m}$ is an uncancelled decomposition of the permutation $\sigma$ in the product of Coxeter generators). Using this fact it is not hard to verify that $a({\varphi}_{q}
)=\frac{q}{(q+1)^2}$.
}\end{exam}
\begin{exam}{\rm
The {\em crystal} \cite{kas} is a set $B$ with the maps
$$\tilde{f}_{i},\tilde{e}_{i}:B\rightarrow B\cup \{0\},$$
where $i\in {\Bbb N}$, for which
$$\tilde{f}_{i}(u)=v \Leftrightarrow \tilde{e}_{i}(v)=u\qquad\mbox{for any}\ u,v\in B,$$
and the functions 
$$\phi_{i}(u) = max\{k\geq 0, \tilde{f}_{i}^{k}(v)\not=0\},$$
$$\varepsilon_{i}(u) = max\{k\geq 0, \tilde{e}_{i}^{k}(v)\not=0\}$$
has finite value for any $u\in B$.
\newline
The morphism $f$ of crystals $B_1$ and $B_2$ is a map of the sets $f:B_1\rightarrow B_2\cup\{0\}$ such that
$$f\tilde{f}_{i}=\tilde{f}_{i}f,\quad f\tilde{e}_{i}=\tilde{e}_{i}f\quad \forall i.$$
The category of crystals will be denoted by $Crystals$. 
\newline
A {\em tensor product} $B_{1}\otimes B_{2}$ of crystals $B_{1}$ and $B_{2}$ is the set $B_{1}\times B_{2}$ with the maps
$$\tilde{f}_{i}(u\otimes v)=
\left\{
\begin{array}{ccc}
& \tilde{f}_{i}(u)\otimes v, & \phi_{i}(u)>\varepsilon_{i}(v) \\
& u\otimes \tilde{f}_{i}(v), & \phi_{i}(u)\leq\varepsilon_{i}(v)
\end{array}\right.$$
$$\tilde{e}_{i}(u\otimes v)=
\left\{
\begin{array}{ccc}
&\tilde{e}_{i}(u)\otimes v, & \phi_{i}(u)\leq \varepsilon_{i}(v) \\
&u\otimes \tilde{e}_{i}(v), & \phi_{i}(u)<\varepsilon_{i}(v)
\end{array}\right.$$
It is easy to see that the tensor product of two crystals is also a cristal and that the category $Crystals$ is monoidal with identical associativity constraint.
\newline
The element $b\in B$ is called {\em highest weight} if
$\tilde{e}_{i}(b)=0\qquad \forall i$. The set of highest weight elements of the crystal $B$ will be denoted by $B^h$.
It follows from definition of the tensor product that
$$(B_{1}\otimes B_{2})^h \subseteq B_{1}^h\otimes B_{2}.$$
\newline
Let us consider the crystal $X=\{x_i ,i\in{\Bbb N}\}$ with the maps
$$\tilde{e}_{i}(x_j )=\delta_{i,j-1} x_{j-1}\quad
\tilde{f}_{i}(x_j )=\delta_{i,j+1} x_{j+1}.$$
It is easy to verify that any connected component of $X^{\otimes n}$ contains only one highest weight element and two components $B_{1}$ and $B_{2}$ are isomorphic iff
$$\phi(b_1 ) = \phi(b_2 )\ \forall i,\ \mbox{where}\ B_{1}^h =\{b_1\},  B_{2}^h =\{b_2\}.$$
It can be proved by the induction that the sequence $\{\phi(b)\}$, where $b\in (X^{\otimes n})^h$, satisfies to the condition $\sum_i i\phi(b) = n$. 
\newline
In other words indecomposible objects of the monoidal subcategory ${\cal T}_0$ in $Crystals$ generated by the object $X$ are parametrized by the partitions. It follows from the results of \cite{nak} that the Littelwood-Richardson coefficients of tensor p
roduct in ${\cal T}_0$ coincides with the standart.
\newline 
The $k$-linear envelope $Crystals_k$ of the category of crystals is a category with the same objects and whose morphisms are (finite) $k$-linear combinations of the morphisms of $Crystals$
$$Hom_{Crystals_k}(B_{1},B_{2}) = \langle Hom_{Crystals}(B_{1},B_{2}) \rangle_k.$$
The category $Crystals_k$ is a semisimple monoidal category.
In particular, the $k$-linear envelope ${\cal T}_{k,0}$ of the subcategory ${\cal T}_0$ in $Crystals$ is a  category ${\cal T}_k$ with another associativity constraint $\varphi_0$.
\newline
It can be verified directly, using the identifications
$$\Lambda^{n}_{0}X = \{x_{i_{1}}\otimes ...\otimes x_{i_{n}}\in X^{\otimes n}, \ i_{1}<...<i_{n}\},$$
that $a(\varphi_0 )=0$.
}\end{exam}
\begin{prop}
Let $k$ is a field of characteristic zero.
The fibres of the map
$$a:Ms({\cal T}_{k})\rightarrow k$$
over $k^*\setminus\{\frac{1}{2(\cos (\alpha )+1)}, \alpha\in{\Bbb Q}^*\}$ consists of one point.
\end{prop}
Proof:
\newline
Let $\psi$ is a associativity constraint for the category ${\cal T}_{k}$.  
\newline
The direct checking shows that the endomorphism $t=q1-(q+1)p \in End_{\cal T}(X^{\otimes 2})$, where $\frac{q}{(q+1)^2}=a$ and $p$ is a projector over $\Lambda^{2}X$, satisfies to the equations
$$(t+1)(t-q)=0,\quad t_{1}t_{2}t_{1}=t_{2}t_{1}t_{2}.$$
Indeed, the first equation follows from the condition $p^2 =p$ and the second from $\ p_{1}p_{2}p_{1}=ap_{1},\quad p_{2}p_{1}p_{2}=ap_{2}$.
\newline
Using the freedom property of the category ${\cal T}_{k,q}$ we can define the monoidal functor 
$$F_{q}:{\cal T}_{k,q}\rightarrow ({\cal T}_{k},\psi),$$
which sends the generator of the Hecke algebra $H_{2}(q)$ to the endomorphism $t$. This functor sends the objects $X$ and $\Lambda_{q}^{2}X$ of ${\cal T}_{k,q}$ to $X$ and $\Lambda^{2}X$ respectively.
\newline
If $a(\psi )=a\in k^*\setminus\{\frac{1}{2(\cos (\alpha )+1)}, \alpha\in{\Bbb Q}^*\}$, then $q$ is not a nontrivial root of unity and the Grothendieck ring $K_{0}({\cal T}_{k,q})$ coincides with $K_{0}({\cal T}_{k})$. The functor $F_{q}$ induces the homo
morphism $K_{0}({\cal T}_{k})\rightarrow K_{0}({\cal T}_{k})$, which preservs $x$ and $\lambda^{2}x$. Hence it is an equivalence by the proposition \ref{hom}.

\end{document}